\documentclass[11pt]{article}
\usepackage[a4paper,right = 22mm, left=22mm,
 top=25mm, bottom = 25mm]{geometry}
\usepackage{graphicx} 
\usepackage{amsmath}
\usepackage{graphicx}
\usepackage{algorithm}
\usepackage{xcolor}
\usepackage{algorithmic}
\usepackage{amsfonts}

\DeclareMathOperator*{\argmin}{arg\,min}                   

\newcommand{\diag}{\mathrm{diag}}

\newcommand{\eps}{\varepsilon}
\renewcommand{\phi}{\mathbf{\varphi}}

\newcommand{\bfT}{\mathbf{T}}

\newcommand{\bfA}{\mathbf{A}}
\newcommand{\bfC}{\mathbf{C}}

\newcommand{\bfU}{\mathbf{U}}

\newcommand{\bfR}{\mathbf{R}}
\newcommand{\bfI}{\mathbf{I}}
\newcommand{\bfD}{\mathbf{D}}

\newcommand{\bfb}{\mathbf{b}}

\newcommand{\bfH}{\mathbf{H}}

\newcommand{\bfx}{\mathbf{x}}

\newcommand{\bfu}{\mathbf{u}}
\newcommand{\bfy}{\mathbf{y}}

\newcommand{\bfW}{\mathbf{W}}
\newcommand{\bfz}{\mathbf{z}}

\newcommand{\bfs}{\mathbf{s}}
\newcommand{\bfv}{\mathbf{v}}
\newcommand{\bfV}{\mathbf{V}}

\newcommand{\bfQ}{\mathbf{Q}}

\newcommand{\bfS}{\mathbf{S}}
\newcommand{\bfZ}{\mathbf{Z}}
\newcommand{\bfSigma}{\mathbf{\Sigma}}
\newcommand{\bfY}{\mathbf{Y}}

\newcommand{\R}{\mathbb{R}}

\newcommand{\trace}{{\mathop{\mathrm{tr}}}}

\newcommand{\bfPsi}{{\boldsymbol{\Psi}}}

\newcommand{\bfepsilon}{{\boldsymbol{\epsilon}}}

\usepackage{caption}
\usepackage{subcaption}
\usepackage{color}
\usepackage{lipsum}
\usepackage{epstopdf}
\usepackage{url} 
\usepackage{authblk}
\newtheorem{proposition}{Proposition}

\numberwithin{equation}{section}
\title{Randomized flexible Krylov methods for $\ell_p$ regularization}
\author[1]{Malena Sabaté Landman}
\author[2]{Yuji Nakatsukasa}
\affil[1]{Mathematical Institute, Oxford, \textit{malena.sabatelandman@maths.ox.ac.uk}}
\affil[2]{Mathematical Institute, Oxford}

\newcommand{\ignore}[1]{}

\date{}
\begin{document}

\maketitle

\begin{abstract}
\noindent The computation of sparse solutions of large-scale linear discrete ill-posed problems remains a computationally demanding task. A powerful framework in this context is the use of iteratively reweighted schemes, which are based on constructing a sequence of quadratic tangent majorants of the $\ell_2$-$\ell_1$ regularization functional (with additional smoothing to ensure differentiability at the origin), and solving them successively.  
Recently, flexible Krylov-Tikhonov methods have been used to partially solve each problem in the sequence efficiently. However, in order to guarantee convergence, the complexity of the algorithm at each iteration increases with respect to more traditional methods. We propose a randomized flexible Krylov method to alleviate the increase of complexity, which leverages the adaptability of the flexible Krylov subspaces with the efficiency of `sketch-and-solve' methods. 
A possible caveat of the mentioned methods is their memory requirements. In this case, one needs to rely instead on inner-outer schemes. In these scenarios, we propose a `sketch-to-precondition' method to speed up the convergence of each of the subproblems in the sequence. 
The performance of these algorithms is shown through a variety of numerical examples.\\

\noindent \textbf{keywords}: Inverse Problems, Randomized Numerical Linear Algebra, Flexible Krylov Methods, Randomized Krylov Methods
\end{abstract}

\section{Introduction}
In many applications in science and engineering, we need to solve large-scale linear ill-posed inverse problems of the form
\begin{equation}\label{eq:linear_pbm}
\bfA \bfx_{\rm true} + \bfepsilon = \bfb, \quad \text{for} \quad \bfA \in\R^{m\times n},
\end{equation}
where $\bfx_{\rm true}$ is the (unknown) true solution, $\bfA$ represents the forward or measurement model, and the measurements $\bfb$ are corrupted with Gaussian white noise $\bfepsilon$. 

Moreover, many problems of the type \eqref{eq:linear_pbm} are typically ill-posed, in the sense that the noise in the measurements is fatally amplified in the reconstruction. In other words, if $\bfA^{\dagger}$ is the Moore-Penrose pseudoinverse of $\bfA$, the approximate solution $\bfA^{\dagger}\bfb$ is far from the true solution $\bfx_{\rm true}$. This is because the singular values of $\bfA$ have a rapid decay and cluster at zero, without a big gap between consecutive values; and noise amplification occurs in the singular subspace associated to the small singular values of $\bfA$. This is true, for example, in problems arising from the discretization of Fredholm integrals of the first kind; which include examples in medical imaging, such as computed tomography (CT); geophysics, such as seismic tomography imaging, or image and signal processing such as deblurring in astronomy or signal deconvolution, see, e.g. \cite{Hansen2010}. This type of problems also arise in linear regressions in statistics, when there is some degree of collinearity in the predictor variables, see e.g. \cite{kutner2004, Callaghan2008IllPosed}.

To obtain a meaningful approximation to the original solution $\bfx_{\rm true}$ in \eqref{eq:linear_pbm} in such ill-posed problems, one needs to use regularization. In other words, one needs to define a new problem whose solution is closely related to the solution $\bfx_{\rm true}$ of \eqref{eq:linear_pbm}, but with better properties: in this case, one that shows better stability with respect to noise in the measurements. Even though there exist a wide variety of regularization methods (including, for example, variational regularization, spectral filtering methods, or iterative regularization \cite{Hansen2010}); these can be characterized by incorporating explicit or implicit information about the expected solution in the problem formulation. The most well-known examples correspond to assuming smoothness in the solution: this is the case of both explicit Tikhonov (or $\ell_2$) regularization  or the use of iterative methods, such as Krylov subspace methods, using early stopping as a regularization mechanism, see, e.g. \cite{Hansen2010}. In particular, there are very efficient ways of combining both approaches, also known as hybrid methods, which also allow for the automatic selection of the regularization parameters, see \cite{Chung2024review} for a recent review. However, these assumptions can be restrictive, and there is a need for efficient algorithms that can be generalized to other regularization terms.

In this paper, we focus on problems whose solutions are assumed to be (numerically) sparse, possibly under a change of basis. This is naturally the case, for example, in imagining applications where the measured objects are physically sparse, for example stars in astronomical imaging. Moreover, this is also used in automatic model (or feature) selection, in statistics, where sometimes the corresponding system matrix is heavily overdetermined \cite{Tibshirani1996Lasso,Sabate2024gmd}. 

It is well known that sparsity in the solution can be promoted using an explicit $\ell_p$ regularization term, for $0 < p \leq 1$, added to the least-squares problem associated with \eqref{eq:linear_pbm}, see \cite{Candes2007Sparsity}. In particular, we solve:
\begin{equation}\label{eq:l2lp}
    \min_{\bfx} \{\|\bfA \bfx -\bfb \|_2^2+\lambda \| \bfPsi\bfx\|_p^p\},
\end{equation}
where $\bfPsi$ is usually a change of basis and $\lambda>0$ is the regularization parameter. Note that we require that $\mathcal{N}(\bfA)\cap\mathcal{N}(\bfPsi)=\{\bf0\}$, where $\mathcal{N}(\cdot)$ denotes the null space of a matrix, so there exists a unique solution of \eqref{eq:l2lp}.
Since for $0<p<1$ the problem \eqref{eq:l2lp} is non-convex, we will mainly focus on the case where $p=1$. Note that this is also known as the LASSO problem in the statistics literature, coined in \cite{Tibshirani1996Lasso}. This corresponds to the convex relaxation of the $\ell_0$-norm regularization problem, which is NP hard \cite{Candes2007Sparsity}.  Moreover, we will consider a smoother version of problem \eqref{eq:l2lp}, which will be defined in Section \ref{sec:background_IRW}.

Problems of the form \eqref{eq:l2lp} have been solved using different schemes. For example, using the following nonlinear optimization methods: fast iterative shrinkage-thresholding algorithm (FISTA) \cite{beck2009fista} and the Sparse Reconstruction by Separable Approximation (SpaRSA) \cite{Wright2008sparsa}. Alternatively, majorization--minimization schemes can be used to approximate the solution of the $\ell_p$ regularized problem by only solving a sequence of appropriate least-square problems. More specifically, this corresponds to an inner-outer loop of iterations, where: (1) a quadratic majorization of the initial problem is updated at each outer loop of iterations so that it is tangent to the current solution, and (2) the solution is updated in each inner loop of iterations to be the approximate solution of this quadratic tangent majorant. This popular framework has also received the name of iteratively reweighted norm (IRN), e.g. \cite{rodriguez2008IRN}, or iteratively reweighted least squares (IRLS),  e.g \cite[4.5.2]{Bjorck1996LS} and has been used in a variety of applications, which might be combined with efficient Krylov subspace methods, see \cite{rodriguez2008IRN, Daubechies2010IRLS, Huang2017MM} just to name a few. More recently, methods that avoid the inner-outer loops of iterations have also been proposed. These are based on building appropriate flexible \cite{Gazzola2014FAT,Chung2019lp,Gazzola2021IRW} or generalized Krylov subspaces \cite{Lanza2015GKS, Huang2017MM} and only partially solving each of the quadratic subproblems in the sequence. 

In the last couple of years, tools from randomized numerical linear algebra have emerged as a key approach in developing fast least-squares solvers for inverse problems, see \cite{meier2022randomized,Sabate2025randCMRH,chung2025randomized}.
In this work, we propose a novel approach to solve inverse problems with $\ell_p$ regularization which incorporate some of these techniques to increase the efficiency of (flexible) Krylov method in this context.

The contributions of this paper can be summarized as follows:
\begin{enumerate}
\item We present a randomized preconditioning approach to solve the  sequence of reweighted least-squares problems appearing when applying an iteratively-reweighted norm scheme for problems with an $\ell_p$ regularization term. Note that this corresponds to a majoritzation-minimization scheme and is therefore based on an inner-outer loop of iterations. This is particularly relevant when the regularization parameter is unknown and each of the subproblems needs to be solved multiple times for different values of the regularization parameter.
\item We present two new randomized flexible Krylov solvers and study their convergence guarantees. The advantage of these methods is that they have a lower complexity than their standard counterparts. Note that, similarly to standard flexible Krylov solvers, these methods are very efficient, and allow for computing the regularization parameter automatically, but have higher memory requirements compared to their non-flexible counterparts.
\end{enumerate}
Theoretical aspects are discussed, and numerical experiments demonstrate the performance of the new methods in comparison to other standard methods.
 
\section{Background on Iteratively Reweighted Norms (IRN)}\label{sec:background_IRW}

Given a function $f:\mathbb{R}^n\rightarrow \mathbb{R}$, a quadratic tangent majorant $Q_k(\bfx)$ of this function at a point $\bfx_{k-1}$ is a quadratic functional that is an upper bound for $f(\bfx)$ over the whole domain; and such that $f(\bfx_{k-1}) = Q_k(\bfx_{k-1})$ and $\nabla f(\bfx_{k-1}) = \nabla  Q_k(\bfx_{k-1})$.  A common way of solving problems of the form \eqref{eq:l2lp} is to consider an inner-outer scheme where a quadratic tangent majorant of \eqref{eq:l2lp} is constructed at the current approximate solution in each outer loop of iterations, and minimized using an iterative method in the inner loop of iterations to update the approximate solution. 

In particular, a common choice of quadratic tangent majorant for \eqref{eq:l2lp} at the point $\bfx_{k-1}$ is 
\begin{gather*}
    Q^{p}_k(\bfx) = {\|\bfA \bfx-\bfb\|}^{2}_{2}+\lambda {\| \bfW^{p}_k \bfPsi \bfx \|}_{2}^{2} + c^p_k \quad \text{where} \quad \bfW^{p}_k = \bfW^{p} ( \bfPsi \bfx_{k-1})\\
 \text{for} \quad \bfW^{p} (\bfz) = \diag( |(\bfz)_1|^{\frac{p-2}{2}},\dots, |(\bfz)_n|^{\frac{p-2}{2}} ),
\end{gather*}
where $(\bfz)_i$ denotes the $i^{th}$ element of the vector $\bfz$,  $c_{k}^{p}=(1-\frac{p}{2}){\| \bfW^{p}_k \bfPsi \bfx_{k-1} \|}_{2}^{2}$ is a constant that does not depend on $\bfx$, and $\lambda$ 
has absorbed any multiplicative constants. Note that this is not the only option of quadratic tangent majorant, see e.g. \cite{mm}. Moreover, to avoid numerical instabilities generated by the lack of smoothness of \eqref{eq:l2lp} on vectors with any 0-valued component, a smoothing parameter $\tau$ is added to the weights in the diagonal of $\bfW^p(\bfz)$:
\begin{equation}\label{eq:smooth_weights}
\bfW^{p, \tau}(\bfz)= \diag(((\bfz)^{2}_1 + \tau^2)^{\frac{p-2}{4}},\dots, ((\bfz)_n^{2}+ \tau^2)^{\frac{p-2}{4}} ),
\end{equation}
so that the following functional is considered instead
\begin{eqnarray}\label{eq:majorant}
    Q^{p,\tau}_k(\bfx) = {\|\bfA \bfx-\bfb\|}^{2}_{2}+\lambda {\| \bfW^{p,\tau}_k \bfPsi \bfx \|}_{2}^{2}+c_k^{\tau}, \, \text{where} \, \bfW^{p, \tau}_k  = \bfW^{p, \tau} (\bfPsi\bfx_{k-1}).
\end{eqnarray}

Where, again, $\lambda$ has absorbed multiplicative constants and $c_{k}^{\tau}$ is a constant that does not depend on $\bfx$. 

A common approach to find approximate solutions to problem \eqref{eq:l2lp}  is using the Majorization-Minimization scheme summarized in Algorithm \ref{alg:MM}. Note that using $Q^{p,\tau}_k(\bfx)$ in Algorithm \ref{alg:MM} leads to minimizing a sequence of quadratic tangent majorants of a smoothed version of \eqref{eq:l2lp}, defined as:
\begin{equation}\label{eq:smoothed}
f^{p,\tau}(\bfx) = \|\bfA \bfx -\bfb \|_2^2+\lambda \| \bfW^{p, \tau}(\bfPsi \bfx) \bfPsi \bfx\|_2^2,
\end{equation}
where the solution-dependent matrix $ \bfW^{p, \tau}(\bfz)$ is defined in \eqref{eq:smooth_weights}. Then, the sequence of solutions produced by Algorithm \ref{alg:MM} will converge to a stationary point of \eqref{eq:smoothed}, and more specifically to its unique minimizer if $p\geq1$. Note that here we have implicitly used the assumption that $\mathcal{N}(\bfA) \cap \mathcal{N}(\bfPsi) = \{{\bf0}\}$, that $\bfW_k^{p, \tau} $ has strictly positive entries and that $\lambda >0$ is fixed. For more details see, e.g. \cite{Lanza2015GKS, Huang2017MM}.

\begin{algorithm}\caption{Quadratic Majorization-Minimization scheme for sparse solutions}\label{alg:MM}
\hspace*{\algorithmicindent} \textbf{Input}: $\bfx_0$, $k_{\max}$, $p$, $\tau$, $\lambda$\\
\hspace*{\algorithmicindent} \textbf{Output}: $\bfx_{k_{\max}}$
\begin{algorithmic}[1]
\FOR {$k=1 \dots k_{\max}$ 
}
\STATE Construct quadratic tangent majorant $Q^{p,\tau}_k(\bfx)$  in \eqref{eq:majorant} of 
\eqref{eq:smoothed} at $\bfx_{k-1}$.
\STATE Update $\displaystyle \bfx_{k} = \argmin_{\bfx}
Q^{p,\tau}_k(\bfx) $ \COMMENT{Possibly involving an inner-loop of iterations}
\ENDFOR
\end{algorithmic}
\end{algorithm}

It is also important to note that $\bfW^{p,\tau}_k$ is invertible by construction, so if we assume that $\bfPsi$ is invertible,  minimizing each of the subproblems involving $Q^{p,\tau}_k(\bfx)$ is equivalent to solving
\begin{eqnarray}\label{eq:reweighted_std}
\bfx_k &=& \bfPsi^{-1} (\bfW_k^{p,\tau})^{-1} \bfs_k,  \quad \text{where} \nonumber \\
\bfs_k &=&\argmin_{\bfs} \{{\|\bfA \bfPsi^{-1} (\bfW_k^{p,\tau})^{-1} \bfs-\bfb\|}^{2}_{2}+\lambda {\|\bfs \|}_{2}^{2}\}. 
\end{eqnarray}
This is true, e.g. when $\bfPsi$ is a change of variables. Alternatively, for example when $\bfPsi$ is a discrete approximation of a differential operator, such as in the case of total variation regularization, other strategies based  on considering $A$-weighted pseudoinverses can be considered and the rest of the algorithm can be used seamlessly, see, e.g. \cite{Sabate2019TV, jimaging7100216}. However, it is outside the scope of this paper to discuss this. 

Even though all derivations follow for $\bfW_k = \bfW^{p,\tau}_k$, from now on we consider only $\bfW_k = \bfW^{1,\tau}_k$ to simplify the notation. Moreover, the regularization parameter $\lambda$, which is typically unknown, will be estimated from the data and can change throughout the iterations. 

\subsection{Regularization parameter choices for Tikhonov regularization}\label{section:reg_param_1}
It is usually the case that the regularization parameter $\lambda$ in \eqref{eq:reweighted_std} is not known a priori. In this case, it is usually estimated from the data using a choice of regularization parameter criterion, see, e.g. \cite{Chung2024review,Gazzola2020param}. Moreover, most regularization parameter choices can be expressed as the minimization of a cost functional involving the solution of \eqref{eq:reweighted_std}, which is now considered to be a function which depends on the parameter, i.e. $\bfx_k(\lambda)$. Therefore, finding a good regularization parameter corresponds to optimizing some functional
\begin{equation}\label{eq:rep_param_optim}
\min_{\lambda >0} P(\bfx_k(\lambda))
\end{equation}
where $\bfx_k(\lambda)=\bfPsi^{-1} (\bfW_k^{p,\tau})^{-1} \bfs_k(\lambda)$ and $\bfs_k(\lambda)$ is the solution of \eqref{eq:reweighted_std} for a given $\lambda$.

One of the most well-known parameter choice criteria is the discrepancy principle, which consists of finding $\lambda_k$ such that
\begin{equation}\label{eq:dp}
     \|\bfA \bfx_k(\lambda_k) - \bfb\|_2 = \tau_{\lambda}\|\bfepsilon\|_2 ,
\end{equation}

where $\bfepsilon$ is the noise present in the measurements \eqref{eq:linear_pbm} and $\tau_{\lambda}>1$ is a safety parameter to avoid over-regularization, see \cite{Chung2024review,Morozov1966OnTS}.
Even though $\|\bfepsilon\|$ is usually not known in practice, we assume that an estimation of the noise level, defined as 
\begin{equation}
    \label{noiselevel}
 \mathrm{nl} = \|\bfepsilon\|_2 / \|\bfb -\bfA \bfx_{\rm true}\|_2,    
\end{equation}
is available, where we usually approximate $ \|\bfb -\bfA \bfx_{\rm true}\|_2 \approx \|\bfb\|_2$.

Alternatively, if the noise level is not known, other regularization parameter choice criteria can be seamlessly used. For example, the generalized cross-validation (GCV) method, which consists on minimizing \eqref{eq:rep_param_optim} where $P$ is the GCV function defined as
\[
G(\lambda) = \frac{\|\bfA \bfx_k(\lambda) - \bfb\|_2}{\trace(\bfI-\bfA_k \bfA_k^{\dagger}(\lambda))},
\]
for 
\[
\bfA_k = \bfA \bfPsi^{-1} \bfW_k^{-1}, \quad
\bfA_k^{\dagger}(\lambda) = \bfA_k (\bfA_k^T \bfA_k + \lambda \bfI)^{-1} \bfA_k^T.
\]

A Majorization-Minimization scheme to approximate the solution of \eqref{eq:l2lp} with a regularization parameter choice can be found in Algorithm \ref{alg:MM_reg_param}. Note that for non-stationary regularization parameters, the convergence of the solutions produced by Algorithm \ref{alg:MM_reg_param} is not well understood.

\begin{algorithm}\caption{Quadratic Majorization-Minimization scheme for sparse solutions with variable regularization parameter $\lambda$}\label{alg:MM_reg_param}
\hspace*{\algorithmicindent} \textbf{Input}: $\bfx_0$, $k_{\max}$, $p$, $\tau$, $\lambda$, stopping criterion. \\
\hspace*{\algorithmicindent} \textbf{Output}: $\bfx_{k_{\max}}$
\begin{algorithmic}[1]
\FOR {$k=1 \dots k_{\max}$}
\STATE Construct quadratic tangent majorant $Q^{p,\tau}_k(\bfx)$ of \eqref{eq:smoothed} at $\bfx_{k-1}$. 
\WHILE{  a stopping criterion is not satisfied} 
\STATE Update $\lambda_{k}$ using a step of a non-linear solver for \eqref{eq:rep_param_optim}.
\STATE Update $\displaystyle \bfx_{k}(\lambda_k)= \argmin{}_{\bfx}\,Q^{p,\tau}_k(\bfx)$. \COMMENT{Possible inner-loop of iterations}
\ENDWHILE
\ENDFOR
\end{algorithmic}
\end{algorithm}

\section{Background on sketching methods}\label{sec:sketching}
A particular technique from randomized numerical linear algebra that has recently gained significant attention in the solution of (possibly regularized) least-squares problems is sketching, see e.g. \cite{Martinsson_Tropp_2020}. This is based on the idea of subspace embedding, i.e., we say that a matrix $\bfS \in \mathbb{R}^{s \times m}$ is a good sketch for a subspace $\mathcal{A}$ if for any $\bfx\in\mathcal{A} \subset \mathbb{R}^{m}$,
\begin{equation}\label{eq:subspace_embedding}
    (1-\eps) \|\bfx\|_2 \leq \|\bfS\bfx\|_2\leq (1+\eps)\|\bfx\|_2,
\end{equation}
where $\eps$ is a distortion parameter. Note that whenever  $s \ll m$, this is a linear dimensionality reduction technique. Fortunately, there are many choices of random matrices which satisfy this property with high probability (for a formal definition, see e.g.~\cite[\S~8.1]{Martinsson_Tropp_2020}), and can be constructed and applied efficiently.  Among these are sparse sketches, SRFT matrices, and leverage-score sampling~\cite[\S~8.3]{Martinsson_Tropp_2020}. 

In this paper, we mainly adopt leverage scores sampling \cite{drineas2012fast,Kolda2022sampling}, which consists of a weighted row selection
as follows. First, compute approximate leverage scores $p_i$ of every row or $\bfA$. Then, draw $s$ samples $\xi_i$ of a multinomial distribution where the probabilities are the normalized approximate leverage scores. Then,
\[
\bfS_{i,j} = \begin{cases}
\sqrt{\frac{\sum_k p_k}{s \, p_j}}, &\text{if} \quad i = \xi_j\\
0 & \text{if} \quad otherwise
\end{cases} \quad i=1,\dots,s\, \,\, j=1,\dots,m.
\]
Given a $k$-dimensional subspace spanned by the columns of an $n\times k$ matrix, 
approximating its leverage scores can be done using (another) sketching and row-norm estimation in $O(nk\log n)$ operations using an SRFT sketch, or $O(nk)$ with a sparse sketch~\cite{drineas2012fast}. Leverage-score sampling has the appeal that (once they are computed or estimated) applying $\bfS$ to a matrix $\bfA$ involves only a row subset of $\bfA$. Here we have an additional reason to use such subsampling-based sketches due to the nature of the problems we are solving, as we describe in Section~\ref{sec:flexible_s2s}, namely equation~\eqref{eq:trick}.

\section{Sketch-to-precondition for IRN schemes}\label{sec:precIRN}
Consider the minimization problem in \eqref{eq:reweighted_std}. First, note that we can write \eqref{eq:reweighted_std} in the augmented least squares form
\begin{equation}\label{eq:reweighted_std_LS}
\bfx_k =  \bfPsi^{-1} \bfW_k^{-1} \bfs_k \quad \text{where} \quad 
\bfs_k =\argmin_{\bfs} \left\|
\begin{bmatrix}
    \bfA \bfPsi^{-1} \bfW_k^{-1} \\
    \sqrt{\lambda} \bfI_n
\end{bmatrix}\bfs - 
\begin{bmatrix}
    \bfb \\ \bf0 
\end{bmatrix}\right\|_2^2.
\end{equation}

Working with the expression \eqref{eq:reweighted_std_LS} is preferable when both the change of variables $\bfPsi$ and the weights $\bfW_k$ are easily invertible, since it reduces instabilities caused by potential large numbers in $\bfW_k$. Note that we have a sequence of problems in standard form, so following from \cite{meier2022randomized}, we can build a good randomized preconditioning for \eqref{eq:reweighted_std_LS} if $\bfA$ is tall and skinny (or fat and short). For simplicity, here we only consider the case where $\bfA\in\mathbb{R}^{m\times n}$ for $m \gg n$, and the other case can be easily generalized from \cite{meier2022randomized}. Consider the following ‘partly exact’ sketch:
\begin{equation*}
\begin{bmatrix}
   \bfS  &  \bf0\\
   \bf0  &  \bfI_n
\end{bmatrix}
\begin{bmatrix}
    \bfA \bfPsi^{-1} \bfW_k^{-1} \\
    \sqrt{\lambda} \bfI_n
\end{bmatrix} = 
\begin{bmatrix}
   \bfS \bfA \bfPsi^{-1} \bfW_k^{-1} \\
    \sqrt{\lambda}  \bfI_n
\end{bmatrix}.
\end{equation*}

Then, a good preconditioner for this system can be built by performing a Cholesky factorization of the Gram matrix of the sketch, 
\begin{equation}\label{eq:chol}
    \bfR_k^T \bfR_k = (\bfS \bfA \bfPsi^{-1} \bfW^{-1}_k)^T (\bfS\bfA \bfPsi^{-1} \bfW_k^{-1}) + \lambda \bfI_n.
\end{equation} 
Note that, if the regularization parameter is $\lambda >0$, we can guarantee that such factorization exists; in finite-precision arithmetic, taking $\lambda> O(u\|\bfS\bfA \bfPsi^{-1} \bfW_k^{-1}\|_2)$ guarantees that Cholesky succeeds. 

Given that the matrix in the right-hand side of \eqref{eq:chol} is of dimension $n \times n$ for small $n$, computing the Cholesky factorization for any choice $\lambda>0$ of the regularization parameter is cheap. Therefore, the ‘partly exact’ sketch can be used to accelerate the convergence of all the solves in the inner-loop of  Algorithm \ref{alg:MM_reg_param}, namely step 5. Moreover, the sketching of $\bfA \bfPsi^{-1}$  can also be re-used for different $\bfW_k^{-1}$ by taking $$\bfY_0 = \bfS \bfA \bfPsi^{-1}, \quad \bfC_0=\bfY_0^T \bfY_0$$ so that $$\bfY = \bfS \bfA \bfPsi^{-1} \bfW^{-1}_k, \quad \bfC=\bfY^T \bfY = \bfW_k^{-1} \bfC_0 \bfW_k^{-1}$$ so then we just need to compute $\text{chol}(\bfC+\lambda \bfI_n)$ for the different weights. This is summarized in Algorithm \ref{alg:sketch-prec}. 

\begin{algorithm}\caption{Randomized preconditioning for IRN schemes}\label{alg:sketch-prec}
\hspace*{\algorithmicindent} \textbf{Input}: $\bfA$, $\bfPsi$, $\bfb$, $\bfx_0$,  $k_{\max}$, $p$, $\tau$, $\lambda$, stopping criterion, sketch $\bfS$.\\
\hspace*{\algorithmicindent} \textbf{Output}: $\bfx_{k_{\max}}$
\begin{algorithmic}[1]
\STATE Construct $\bfY_0=\bfS\bfA\bfPsi^{-1}$, and  $\bfC_0 = \bfY_0^{T} \bfY_0$.
\FOR {$k=1\dots k_{\max}$}
\STATE Construct $\bfW_k =  \bfW^{p, \tau}(\bfPsi \bfx_{k-1})$  as defined in \eqref{eq:smooth_weights}.
\WHILE{  a stopping criterion is not satisfied} 
\STATE Compute $\bfR_k = \text{chol}(\bfW_k^{-1} \bfC_0\bfW_k^{-1} + \lambda_k \bfI_n)$.
\STATE Update $\lambda_k$ using a step of a non-linear solver for \eqref{eq:rep_param_optim}.
\STATE Update $\bfx_{k}(\lambda_k)$  by solving \eqref{eq:reweighted_std_LS} with right-preconditioning $\bfR_k^{-1}$.
\ENDWHILE

\ENDFOR
\end{algorithmic}
\end{algorithm}

\section{Randomized flexible Krylov subspaces}
If there is enough memory to store a moderate amount of basis vectors for the Krylov subspace, one can consider the weights in expression \eqref{eq:reweighted_std} to be iteration-dependent preconditioners, also called priorconditioners \cite{calvetti2005priorconditioners}, and use flexible Krylov subspaces, see e.g. \cite{Chung2019lp, Gazzola2021IRW} and references therein. These methods avoid the inner-outer loops of iterations in traditional majorization-minimization schemes by building a single solution subspace of increasing dimension, spanned by the columns of a matrix that we call $\bfZ_k$, in which the sequence of regularized subproblems of the form \eqref{eq:reweighted_std} are (approximately) solved. 

The key distinction between a classical Krylov subspace and a flexible Krylov subspace is that in the flexible case, the preconditioning matrix (in this case $\bfW_k^{-1}$) can change at each iteration. To construct such spaces, one can use either the flexible Golub-Kahan (FGK) or the flexible Arnoldi (FA) factorizations:
\begin{eqnarray*}\label{flex_gk_A}
    &\text{(FGK)}&\quad  \bfA \bfPsi^{-1} \bfZ^{\text{\tiny FGK}}_k = \bfU^{\text{\tiny FGK}}_{k+1} \bfH^{\text{\tiny FGK}}_{k+1,k}, 
    \quad 
    \bfA^\top \bfU^{\text{\tiny FGK}}_{k+1} = \bfV^{\text{\tiny FGK}}_{k+1} \bfT^{\text{\tiny FGK}}_{k+1,k},\\
    &\text{(FA)}&\quad  \bfA \bfPsi^{-1} \bfZ^{\text{\tiny FA}}_k = \bfV^{\text{\tiny FA}}_{k+1} \bfH^{\text{\tiny FA}}_{k+1,k},
\end{eqnarray*}  
which we write in unified form as  
\begin{equation}\label{flex_exp}
    \bfA \bfPsi^{-1} \bfZ_k = \bfU_{k+1} \bfH_{k+1,k}, \quad  \text{where} \quad   \bfZ_k = \big[ \bfW_1^{-1} \bfv_1, \dots, \bfW_k^{-1} \bfv_k \big].
\end{equation}  

Note that, for flexible Arnoldi $\bfU_{k+1}=\bfV_{k+1}$. Unlike in the classical case, the columns of $\bfZ_k$, which span the flexible Krylov solution subspace,  are not orthonormal because of the changing preconditioners. In contrast, $\bfU_{k+1}$ and $\bfV_{k+1}$ have orthonormal columns and $\bfH_{k+1,k}$ is upper Hessenberg. See, e.g., \cite{Chung2019lp,Gazzola2021IRW}, for more information and explicit algorithms for the solution subspace updates.

The specific choice of flexible factorization and the corresponding optimality conditions determine which flexible Krylov method is applied and, in this work, we focus on minimal residual methods. Flexible Krylov methods have been used to find the solution of linear systems of the form \eqref{eq:linear_pbm} in different ways that can be summarized in the following general minimization problem:
\begin{equation}\label{eq:flex_opt}
    \min_{\bfy} \{\|\bfA \bfPsi^{-1} \bfZ_k \bfy - \bfb\|_2^2 + \lambda R_k(\bfy) \},
\end{equation}
where $\lambda R_k(\bfy)$ is a regularization term that needs to be specified. Moreover, since the iteration-dependent preconditioning $\bfW_k^{-1}$ cannot be easily factored out of the matrix containing the basis vectors, i.e. $\bfZ_k \neq \bfW_k^{-1} \bfV_k$, we find a preconditioned solution directly. This means, if the weights in the preconditioner were fixed, $\bfW_k^{-1} = \bfW^{-1}$, then following the notation in \eqref{eq:reweighted_std_LS}, $\bfZ_k \bfy_k = \bfW_k^{-1}  \bfs_k$, for $\bfs_k$ in \eqref{eq:reweighted_std_LS}, so we are implicitly working with a change of variables. 

Problem \eqref{eq:flex_opt} can be easily solved using the relationship in \eqref{flex_exp} to efficiently compute the projection, and this has the advantage of allowing the regularization parameter to be chosen on-the-fly. The difference between the different approaches is the following. The first approach consists in setting $R_k(\bfy)=0$ and regularizing by early stopping of the iterations, as \cite{Gazzola2014FAT}. Let us assume that the matrix-vector multiplication can be implemented in $O(n)$ operations. For the GMRES-based method, the complexity is $O(n k^2)$ (since, compared to GMRES, it only requires an extra matrix vector product with a diagonal matrix per iteration; here we are assuming matrix-vector multiplication with $\bfA$ is efficient, e.g. $O(n)$). However, without early stopping, the approximated solution computed with this solver semiconverges to the solution $\bfA^{\dagger}\bfb$. Secondly, following a `project-first-then-regularize' approach, also known as a hybrid method, one can consider $R_k(\bfy)=\|\bfy\|_2^2$, see, e.g. \cite{Chung2019lp}. For the GMRES-based algorithm, this approach keeps the complexity at $O(n k^2)$ (since, compared to FGMRES, it only requires regularizing the ($k$+$1$) $\times$ $k$  projected problem) and it has been empirically observed that early stopping of the iterations provides a `sparse solution'. However, without stopping the iterations, the solution will `semi-converge' to the Tikhonov-regularized least squares solution. Lastly, one can use an iteratively-reweighted scheme, where $R_k(\bfy)=\|\bfW_k \bfZ_k \bfy\|_2^2$, as in \cite{Gazzola2021IRW}. This method has theoretical convergence guarantees, but requires the computation of the QR factorization of the tall and skinny matrix $\bfW_k \bfZ_k$ at each iteration. Moreover, this cannot be cheaply updated at each iteration as $\bfZ_k$ increases dimensions, even if only a new column is appended to it, i.e. $\bfZ_k =[\bfZ_{k-1},\bfz_k]$, because $\bfW_k$ is also changing. This increases the complexity of the algorithm to $O(n k^2)$ + $O(n k^3)$, where the second term relates to a `fresh' QR factorization required at each iteration. 

\subsection{Sketch-and-solve for flexible Krylov subspaces}\label{sec:flexible_s2s}
To reduce the computational cost of flexible Krylov subspace methods, one can resort to randomization. First, we consider a `sketch-and-solve' approach, where, at each iteration, we truncate the expensive re-orthogonalization in the flexible Arnoldi or Golub-Kahan processes, and we minimize a sketched version of \eqref{eq:flex_opt}, as is detailed later in the section.

The complexity of sketching a matrix $n \times k$ is $O(n\,k \, \log k)$ with an SRFT sketch or $O(nk)$ with a sparse sketch, so we will consider the latter. Moreover, assume that $O(n)=O(m)$. Then, using sketching in the regularization term already reduces the complexity of the IRW flexible methods from $O(n k^2)+O(n k^3)$ to $O(n k^2)+O(n k) +O(k^4)$. Moreover, we extend the sketch-and-solve approach to the fit-to-data term to further reduce the complexity to $O(nk)$ + $O(k^4)$ + $O(k^3)$, avoiding the (full) re-orthogonalization of the columns of $\bfU_{k+1}$ and $\bfV_{k}$ in \eqref{flex_exp}. 

Sketch-and-solve flexible Krylov subspace methods consist of the following. Assume that, at each iteration $k$, we have performed $k$ steps of the $\ell$-truncated flexible Arnoldi or the $\ell$-truncated flexible GK processes. This means, we replace the full orthogonalization of $\bfU_{k+1}$ and $\bfV_k$ in the FA and FGK factorizations, summarized in \eqref{flex_exp} and defined above, by only orthogonalizing each vector $\bfu_{k+1}$ ($\bfv_k$) against the last $\ell$ columns of $\bfU_{k}$ ($\bfV_{k-1}$). Note that, in this case, we obtain a relation which is formally equivalent to \eqref{flex_exp},
\begin{equation*}
    (\ell\text{-FGK} / \ell\text{-FA}) \quad     \bfA \bfPsi^{-1} \bar\bfZ_k = \bar\bfU_{k+1} \bar\bfH_{k+1,k}, \quad  \text{where} \quad   \bar\bfZ_k = \big[ \bfW_1^{-1} \bar\bfv_1, \dots, \bfW_k^{-1} \bar\bfv_k \big],
\end{equation*}
so that the corresponding solution space is spanned by the columns of $\bar{\bfZ}_k$, but where $\bar\bfU_{k+1}$ and $\bar\bfV_k$ do not have orthonormal columns. Throughout the algorithm, we are implicitly (and iteratively) computing the matrix on the left-hand-side of ($\ell\text{-FGK} / \ell\text{-FA}$):
\begin{eqnarray}\label{truncatedGKB}
   \bfA\, \bfPsi^{-1} \bar{\bfZ}_k,
\end{eqnarray}
and a sketch of \eqref{truncatedGKB} can then be used to approximate the residual and compute an approximated minimal residual norm solution. It is important to note that the space spanned by the columns of $\bar{\bfZ}_k$ (coming from a flexible truncated process) and $\bfZ_k$ (coming from a flexible process with full orthogonalization), might not be the same for iteration dependent preconditioning $\bfW_k^{-1}$, since the weights will be applied to different vectors. However, we expect either basis to promote the desired features in our solution. It is important to remark that, since we need to compute the solution at each iteration, we need to store all basis vectors at every iteration. However, since we are only considering $\ell$-truncated methods, we only need $\ell$ vectors to be stored in working precision (plus a sketched version of the vectors, as is explained later), and we can store the first $k-\ell$ vectors in very low precision to obtain an approximation of the solutions. 

Once we have computed a non-orthogonal basis for the coefficients, $\bfPsi^{-1} \bar{\bfZ}_k$, and our residuals, $\bfA \bfPsi^{-1} \bar{\bfZ}_k$, as in \eqref{truncatedGKB},
we consider a `sketch-and-solve' process to approximate a solution of the following problem
\begin{equation}\label{eq:min_sketch_about_to_project}
    \min_{\bfy} \{\|\bfA \bfPsi^{-1} \bar{\bfZ}_k \bfy - \bfb\|_2^2 + \lambda_k \|\bfW_k \bar{\bfZ}_k\bfy\|_2^2 \},
\end{equation}
so we consider the QR factorization of the short sketched matrices 
\begin{eqnarray} \label{eq:R1R2}
\bfS_1 \bfA \bfPsi^{-1} \bar{\bfZ}_k &=&  \bfQ_1 \bfR_1 \nonumber \\
\bfS_2 \bfW_k \bar{\bfZ}_k &=&   \bfQ_2 \bfR_2, 
\end{eqnarray}
where $\bfS_1$, $\bfS_2$ are sketching matrices and the updates for weights can be written as $\bfW_{k} = \bfW^{p,\tau}(\bfPsi \bfx_{k-1})$, defined in \eqref{eq:smooth_weights}. 

Now, at each iteration $k$ of the proposed sketch-and-solve flexible Krylov method, an approximate solution to \eqref{eq:l2lp} is sought by sketching-to-project equation \eqref{eq:min_sketch_about_to_project}:
\begin{equation}\label{sketch_and_solve}
    \bfx_k = \bfPsi^{-1} \bar{\bfZ}_k \bfy_k \quad \text{for} \quad \bfy_k=\argmin_{\bfy} \{\|\bfR_1 \bfy - \bfQ_1^T\bfS_1\bfb\|_2^2 + \lambda_k \|\bfR_2\bfy\|_2^2 \}.
\end{equation}
Depending on the choice of the basis vectors, determined by the use of Arnoldi or Golub-Kahan based methods, we call the methods associated to the minimization \eqref{sketch_and_solve} sketch-and-solve iteratively-reweighted flexible GMRES (S\&S-IRW-FGMRES) or LSQR (S\&S-IRW-FLSQR). A summary of the general scheme for both methods can be found in Algorithm \ref{alg:sketch-flex}.

\begin{algorithm}\caption{Sketch-and-solve  flexible Krylov methods}\label{alg:sketch-flex}
\hspace*{\algorithmicindent} \textbf{Input}: $\bfA$,$\bfPsi$,$\bfb$,$\bfx_0$, $k_{\max}$, $p$, $\tau$, $\lambda$, stopping criterion, sketches $\bfS_1$ and $\bfS_2$.\\
\hspace*{\algorithmicindent} \textbf{Output}: $\bfx_{k_{\max}}$
\begin{algorithmic}[1]
\FOR {$k=1\dots k_{\max}$}
\STATE Construct $\bfW_k = \bfW^{p,\tau}(\bfPsi \bfx_{k-1})$ as defined in \eqref{eq:smooth_weights}.
\STATE Increase your search space by computing $\bar{\bfz}_{k}$.
\STATE Update $\bfR_1$ and $\bfR_2$ in \eqref{eq:R1R2}.
\STATE Compute $\lambda_{k}$ using the DP or GCV in the sketched system \eqref{sketch_and_solve}.
\STATE Update $\bfx_{k}$ by solving \eqref{sketch_and_solve} for given $\lambda_k$.
\ENDFOR
\end{algorithmic}
\end{algorithm}

\begin{proposition}
Consider the sketch-and-solve flexible LSQR or GMRES algorithms approximating the minimizer of $f^{p,\tau}(\bfx)$ in \eqref{eq:smoothed} after $k \ll \min(m,n)$ iterations, and where no break-down has happened. Assume, without loss of generality, that $\bfPsi=\bfI_n$. Moreover, let the regularization parameter $\lambda$ be fixed, as well as $p \in (0,2]$ and $\tau$.  Last, assume that we have two sketches $\bfS_1$ and $\bfS_2$ that fulfill the subspace embedding property \eqref{eq:subspace_embedding} for the range of $[\bfA \bar\bfZ_k, \bfb]$ and $\bfW_k\bar\bfZ_k$ respectively, where the maximum between their distortion parameters is $0\leq\eps_k<1$. 

Under these assumptions, given a distortion factor $\eps_k$ induced by the sketching, if 
\begin{equation}\label{eq:mono_cond}
\frac{\hat Q^{p,\tau}_k(\bfx_{k-1}) -  \hat Q^{p,\tau}_k(\bfx_{k})}{\hat Q^{p,\tau}_k(\bfx_{k})} \geq \frac{2\eps_k}{1-\eps_k}
\end{equation}
where $\hat Q^{p,\tau}_k(\bfx)$ is the sketched functional then,
$$f^{p,\tau}(\bfx_{k-1}) \geq f^{p,\tau}(\bfx_{k}) \geq 0.$$  
\end{proposition}

This result shows that, for sufficient relative decrease in the sketched functional, the value of the original functional evaluated at the solution $\bfx_k$  decreases monotonically and is bounded from below by zero. 

\textbf{Proof.} Consider the quadratic tangent majorant  $Q^{p,\tau}_k(\bfx)$ of $f^{p,\tau}(\bfx)$ at $\bfx_{k-1},$ as defined in \eqref{eq:majorant} and define the sketched functional
\[
    \hat Q^{p,\tau}_k(\bfx)=\| \bfS_1 (\bfA \bfx - \bfb)\|_2^2 + \lambda \|\bfS_2 \bfW_k \bfx\|_2^2.
\]
By assumption, if $\eps_k$ is the maximum between the  distortion parameters associated to the sketches $\bfS_1$ and $\bfS_2$, then 
\begin{equation}\label{eq:sketch_bounds}
(1-\eps_k) \,(Q^{p,\tau}_k(\bfx)- c^p_k) \leq \hat Q^{p,\tau}_k(\bfx) \leq (1+\eps_k) \,(Q^{p,\tau}_k(\bfx)- c^p_k) \quad  \forall \bfx\in\mathcal{R}(\bar\bfZ_{k})
\end{equation} 
where the distortion factor $\eps_k$ depends on the sketch. Moreover, for (fixed) given sketches $\bfS_1$ and $\bfS_2$, 
\begin{eqnarray*}
    \hat Q^{p,\tau}_k(\bfx_{k}) = \min_{\bfx \in \mathcal{R}(\bar \bfZ_{k})}\hat Q^{p,\tau}_k(\bfx) \leq \hat Q^{p,\tau}_k(\bfx_{k-1}),
\end{eqnarray*}
since $\mathcal{R}(\bar\bfZ_{k-1}) \subseteq \mathcal{R}(\bar\bfZ_{k})$, so we can define 
$\Delta =  \hat Q^{p,\tau}_k(\bfx_{k-1}) -  \hat Q^{p,\tau}_k(\bfx_{k}) \geq 0. $
Now, using \eqref{eq:sketch_bounds} we can consider
\begin{eqnarray*}
 Q^{p,\tau}_k(\bfx_{k-1}) -  Q^{p,\tau}_k(\bfx_{k}) =  (Q^{p,\tau}_k(\bfx_{k-1}) - c_k^p) - (Q^{p,\tau}_k(\bfx_{k})- c_k^p) \\
  \geq \frac{1}{1+\eps_k}\hat Q^{p,\tau}_k(\bfx_{k-1}) - \frac{1}{1-\eps_k}\hat Q^{p,\tau}_k(\bfx_{k}) \\
  \geq \frac{ (1-\eps_k) ( \Delta+ \hat Q^{p,\tau}_k(\bfx_{k}) )
- (1+\eps_k)\hat Q^{p,\tau}_k(\bfx_{k})}{(1+\eps_k)(1-\eps_k)},
\end{eqnarray*}
and this is non-negative if 
\[(1-\eps_k) ( \Delta+ \hat Q^{p,\tau}_k(\bfx_{k}) )
- (1+\eps_k)\hat Q^{p,\tau}_k(\bfx_{k}) \geq 0 \]
which holds if the relative step is greater than the following tolerance:
\[\frac{\Delta}{\hat Q^{p,\tau}_k(\bfx_{k})} \geq \frac{2\eps_k}{1-\eps_k}.\]
where we have used that $\hat Q^{p,\tau}_k(\bfx_{k}) > 0 $ (if this is 0, we have found a minimum, since the function is non-negative).

Using this, we have that
\[f^{p,\tau}(\bfx_{k}) \leq Q^{p,\tau}_k(\bfx_{k}) \leq Q^{p,\tau}_k(\bfx_{k-1}) = f^{p,\tau}(\bfx_{k-1}), 
\]
where the first inequality holds since $ Q^{p,\tau}_k(\bfx)$ is a majorant of $f^{p,\tau}(\bfx)$ and the last equality holds given the definition of $ Q^{p,\tau}_k(\bfx)$. \hfill $\square$

\textbf{Remark.} 
In the case of sketches where the distortion can be empirically estimated, or it is known analytically with high probability, one can use information from the projected problem to assess if the original full-dimensional function is decreasing. Moreover, this could be used to modify the algorithm when decrease cannot be guaranteed by considering a bigger sketch, or triggering a restart.

We should also note that $\eps_k$ is often not very small, say $0.5$, so the assumption~\eqref{eq:mono_cond} is often stringent. We revisit this in Proposition~\ref{prop2} when the least-squares problems are solved to full accuracy. 

\subsubsection{Efficient sketching} In this subsection, we specify a very efficient sketching procedure that can be updated throughout the iterations. First, assuming that we have a maximum number of iterations $k_{\max}$, we set $\bfS_1$ and $\bfS_2$ to have $O(k_{\max})$ rows.  For the first line of \eqref{eq:R1R2}, note that only one column is appended to $\bar{\bfZ}_k$ at each iteration $k$, so that only one column has to be appended to the sketched matrix, and the QR factorization can be correspondingly cheaply updated. The second line of \eqref{eq:R1R2} is less straightforward, since the matrix $\bfW_k$ changes at each iteration. However, we know that this is a diagonal matrix so, in particular, for $\bfS_2 \in \mathbb{R}^{s \times n}$ we have that
\begin{equation}\label{eq:trick}
\bfS_2 \bfW_k = \bar{\bfW}_k \bfS_2, \quad \text{provided that}  \quad  \bar{\bfW}_k = \bfS_2 \bfW_k \bfS^{\dagger}_2 =  \bfS_2 \bfW_k \bfS_2^T  (\bfS_2 \bfS_2^T )^{-1} ,
\end{equation}
for any $\bfS_2$ with full row rank. If we consider $\bfS_2$ to be a sparse sketching matrix with this property, then 
\[
\bfS_2 \bfW_k \bar{\bfZ}_k = \bar{\bfW}_k \bfS_2 \bar{\bfZ}_k  =  \bfQ_2 \bfR_2, 
\]
so we can re-use and only update the sketch $\bfS_2 \bar{\bfZ}_k$. Therefore, if we can keep the complexity of computing $\bar{\bfW}_k$ to order $O(n)$ at each iteration, then the complexity of the sketching is $O(kn)$. The most well known sketch that naturally satisfies this property~\eqref{eq:trick}  is row selection, where the rows can be selected, or sampled, using different techniques. In this paper, we use leverage score sampling as explained in Section \ref{sec:sketching}. 

\subsubsection{Regularization parameter choices}
Note that, since \eqref{sketch_and_solve} is a problem of small dimension, this allows for the regularization parameter $\lambda_k$ to be chosen on-the-fly using a regularization parameter choice method. For example, we can use the discrepancy principle as defined in \eqref{eq:dp}, where we now use sketching to approximate the residual norm. In particular, consider the solution for the general Tikhonov problem \eqref{sketch_and_solve}, i.e.
\begin{equation}\label{eq:genTik}
\bfy_k(\lambda) = (\bfR_1^{T}\bfR_1+\lambda \bfR_2^T\bfR_2)^{-1} \bfR_1^{T} {\bfQ}_1^{T}\bfS_1\bfb =  \bfR_{\lambda}^{\dagger} {\bfQ}_1^{T}\bfS_1\bfb .  
\end{equation}

Then, $\lambda_k$ is chosen to satisfy:
\[
 \|\bfQ_1\bfR_1 \bfy_k(\lambda_k) - \bfS_1\bfb\|_2  = \tau_{\lambda} \mathrm{nl} \|\bfb\|_2 ,
\]   
(where $\mathrm{nl}$ is the noise level~\eqref{noiselevel}) or $\lambda_k=0$ if the above equation does not have a solution – this makes sense, as the minimal projected norm on the left-hand-side is bigger than the right-hand-side, so $\lambda_k=0$ minimizes the difference between the two. 

As explained in Section \ref{section:reg_param_1}, if the noise level is not known, one can use other criteria. For example, considering $\bfR_{\lambda}^{\dagger}$ as defined in \eqref{eq:genTik}, the GCV function for the sketched problem is defined as 
\begin{equation}\label{eq:proj_GCV}
G_k(\lambda) = \frac{k\|\bfR_1 \bfy_k(\lambda)-{\bfQ}_1^T \bfS_1 \bfb\|_2^2}{(\trace(\bfI-\bfR_1 \bfR_{\lambda}^{\dagger}))^2}=
\frac{k\|(\bfI-\bfR_1 \bfR_{\lambda}^{\dagger}) {\bfQ}_1^{T}\bfS_1\bfb \|_2^2}{(\trace(\bfI-\bfR_1 \bfR_{\lambda}^{\dagger}))^2},
\end{equation}
which can be easily computed using the GSVD of the small matrices $\bfR_1$ and $\bfR_2$. Then, at each iteration $k$, we take the regularization parameter $\lambda_k$ which minimizes the projected GCV function \eqref{eq:proj_GCV}. However, this regularization parameter might produce overly smooth solutions, so we use the weighted GVC (WGCV) function instead:
\[
G^{\omega}_k(\lambda) = \frac{k\|\bfR_1 \bfy_k(\lambda)-{\bfQ}_1^T \bfS_1 \bfb\|_2^2}{(\trace(\bfI-\omega\bfR_1 \bfR_{\lambda}^{\dagger}))^2}=
\frac{k\|(\bfI-\bfR_1 \bfR_{\lambda}^{\dagger}) {\bfQ}_1^{T}\bfS_1\bfb \|_2^2}{(\trace(\bfI- \omega \bfR_1 \bfR_{\lambda}^{\dagger}))^2},
\]
where $\omega =(k+1) / s$, with $s$ being the size of the sketch, as proposed in \cite{renaut2017hybrid}. Then, at each iteration, we use
\[
\lambda_k = \arg\min_{\lambda} G^{\omega}_k(\lambda).
\]

Note that these strategies are suitable to find appropriate regularization parameters for each of the projected problems. However, other approaches based on approximations of the full-dimensional regularization parameter choice criteria can also be considered, see e.g. \cite{Chung2024review,Gazzola2020param}.

Last, following the convention for hybrid flexible methods explained after equation \eqref{eq:flex_opt}, in the case where $\lambda_k=0$ for all $k$, we call these methods sketch-and-solve flexible GMRES (S\&S-FGMRES) or LSQR (S\&S-FLSQR).

\subsection{Sketch-to-precondition for flexible Krylov subspaces} \label{sec:flexible_s2p}
Lastly, we can solve the projected problems defined in \eqref{eq:min_sketch_about_to_project} using a sketch-to-precondition approach. Note that these are all tall-and-skinny least-squares problems, for which we can use a preconditioning that is easily adaptable to different regularization parameters in the same spirit of Section \ref{sec:precIRN}. More specifically, we use the same sketching as in \eqref{eq:R1R2}, and construct (iteratively)
\begin{equation}\label{eqn:CandD}
\bfC_k\hspace{-2pt}=\hspace{-2pt}(\bfS_1 \bfA \bfPsi^{-1} \bar{\bfZ}_k )^T (\bfS_1 \bfA \bfPsi^{-1} \bar{\bfZ}_k )\hspace{-2pt}=\hspace{-2pt} \bfY_k^T\bfY_k, \, \hspace{8pt}\bfD_k \hspace{-2pt}= \hspace{-2pt}(\bfS_2 \bfW_k \bar{\bfZ}_k  )^T (\bfS_2 \bfW_k \bar{\bfZ}_k) ,    
\end{equation}
which can then be used to update the randomized preconditioner $\bfR_k$ for the projected problem \eqref{eq:min_sketch_about_to_project} as the Cholesky factor of a small $k\times k$ matrix:
\begin{equation}\label{prec_flex}
    \bfR_k^T\bfR_k=\bfC_k+\lambda \bfD_k.
\end{equation}
Note that, however, in this case we need to deal with a problem that is not in standard form. The general structure of these type of methods, which we call S2P-IRW-FGMRES or S2P-IRW-FLSQR depending on the solution subspace they are inherently constructing, is summarized in Algorithm~\ref{alg:sketch-prec-flex}.

Note that sketch-to-precondition methods for standard Krylov solvers are equivalent to their non-preconditioned counterparts. This is because they preserve the same optimality conditions, and the truncated recurrences change the basis vectors but preserve the same solution space. However, this is not true for flexible Krylov methods, i.e. S2P-IRW-FGMRES and S2P-IRW-FLSQR are not equivalent to IRW-FGMRES and IRW-FLSQR, respectively. As already mentioned in the previous subsection, this is the effect of the truncation in the flexible Arnoldi or flexible Golub-Kahan factorizations, and not of the sketch-to-precondition itself. In fact, both methods assume the same optimality conditions, but the combination of the iteration-dependent preconditioning and the truncated flexible Arnoldi or Golub-Kahan methods generate different search spaces for the solution. 

Compared to S\&S-IRW methods, S2P-IRW methods are theoretically more robust, since the optimality conditions allow for stronger monotonicity guarantees. In particular, these are based on the subsequent projections of a sequence of quadratic tangent majorants of the original (smoothed functional) onto a subspace of increasing dimension, in the same spirit as other flexible Krylov methods for IRN schemes, see e.g. \cite{Gazzola2021IRW}. Therefore, we can directly use the theory in e.g. \cite{Gazzola2021IRW} to guarantee a descent step at each iteration.

\begin{algorithm}\caption{Sketch-to-precondition flexible Krylov methods}\label{alg:sketch-prec-flex}
\hspace*{\algorithmicindent} \textbf{Input}: $\bfA$,$\bfb$,$\bfPsi$,$\bfx_0$, $k_{\max}$, $p$, $\tau$, $\lambda$\\
\hspace*{\algorithmicindent} \textbf{Output}: $\bfx_{k_{\max}}$
\begin{algorithmic}[1]
\FOR {$k=1\dots k_{\max}$}
\STATE Construct $\bfW_k = \bfW^{p,\tau}(\bfPsi \bfx_{k-1})$ as defined in \eqref{eq:smooth_weights}.
\IF {$k < \min(m,n)$}
\STATE Increase your search-space by computing $\bar{\bfz}_{k}$, 
\STATE compute $\bar{\bfy}_{k} = \bfS \bfA\bfPsi^{-1}\bar{\bfz}_{k}$,
\ELSE 
\STATE consider $\bar \bfZ_k=\bfI_n$.
\ENDIF
\STATE Update $\bfY_k$ by appending $\bar{\bfy}_{k}$, i.e., $\bfY_k\leftarrow [\bfY_k,\, \bar{\bfy}_{k}]$, and update $\bfC_k = \bfY_k^{T} \bfY_k$.
\STATE Update $\bfD_k$ as in~\eqref{eqn:CandD} by recomputing or using \eqref{eq:trick}.
\STATE Compute $\lambda_{k}$, possibly using an inner-loop of iterations, and updating the preconditioner $\bfR_k=\text{chol}(\bfC_k+\lambda \bfD_k)$ in \eqref{prec_flex} for \eqref{eq:min_sketch_about_to_project} for each $\lambda$.
\ENDFOR
\end{algorithmic}
\end{algorithm}

\begin{proposition}\label{prop2}
Consider the sketch-to-precondition iteratively-reweighted flexible LSQR (S2P-IRW-FLSQR) or GMRES (S2P-IRW-FGMRES) algorithms approximating the minimizer of $f^{p,\tau}(\bfx)$ in \eqref{eq:smoothed}. Then, for fixed $\lambda$, $p \in (0,2]$ and $\tau$, the value of the original functional evaluated at the solution $\bfx_k$  decreases monotonically and is bounded from below by zero. In other words,
$$f^{p,\tau}(\bfx_{k-1}) \geq f^{p,\tau}(\bfx_{k}) \geq 0.$$ 
\end{proposition}
\textit{Proof} Consider the quadratic tangent majorant  $Q^{p,\tau}_k(\bfx)$ of $f^{p,\tau}(\bfx)$ at $\bfx_{k-1},$ as defined in \eqref{eq:majorant}. By definition, we have that
\[
f^{p,\tau}(\bfx_{k}) \leq Q^{p,\tau}_k(\bfx_{k}) \quad \text{and} \quad  Q^{p,\tau}_k(\bfx_{k-1}) = f^{p,\tau}(\bfx_{k-1})
\]
and, given the definition of $\bfx_k$,
\[
Q^{p,\tau}_k(\bfx_{k}) = \min_{\bfx\in \mathcal{R}(\bar\bfZ_{k})} Q^{p,\tau}_k(\bfx) \leq Q^{p,\tau}_k(\bfx_{k-1}).
\]
Here, the inequality holds since $\mathcal{R}(\bar \bfZ_{k-1}) \subseteq \mathcal{R}(\bar \bfZ_{k})$, where the equality only holds after breakdown of the truncated flexible Arnoldi or the truncated flexible GK, or after $k \geq \max(m,n)$.
Putting everything together, we have that
$$f^{p,\tau}(\bfx_{k-1}) \geq f^{p,\tau}(\bfx_{k}),$$ 
and by definition, $f^{p,\tau}(\bfx_{k}) \geq 0$, so that completes the proof.\hfill $\square$ \\

\textbf{Remark.} Note that, similarly to the iteratively reweighted flexible Krylov methods, we have extended the solver so that the weights are still updated at each iteration after the solution space is $\mathbb{R}^n$. After that point, the method becomes the classic IRN, and therefore its convergence theory follows. That is, the method will converge to a stationary point of \eqref{eq:smoothed} and, if $p\geq 1$, it will converge to the unique solution of \eqref{eq:smoothed} \cite{Gazzola2021IRW}.

\section{Numerical experiments} In this section, we report three numerical examples highlighting the performance of the methods introduced in this paper; each of them focusing on a different aspect.

The first experiment considers a simulated subset‐selection problem, illustrating the effectiveness of the sketch-to-precondition method for problems involving a tall and skinny matrix. In this example, we compare methods based on classic Krylov solvers, and study the effect of including the above-mentioned preconditioning. The second experiment corresponds to a deblurring problem with a numerically sparse solution, representing a starry sky, where we compare methods based on the (flexible) Arnoldi process with a fixed regularization parameter. In this example, we focus on the convergence results in Section \ref{sec:flexible_s2s}, as well as a comparison with other hybrid classic and flexible Krylov subspace methods, and other standard solvers for problems with $\ell_1$ regularization.  Last, a simulated computed tomography (CT) example is shown to compare the methods based on the (flexible) Golub-Kahan algorithm. In this example, we aim to illustrate the properties of the different subspace basis vectors and their inherent regularization properties, as well as the benefits of using explicit regularization and the effect of  the regularization parameter choice.

In the following, we summarize the methods compared in this section, along with their chosen abbreviations and their corresponding citations.  
First, we recall the three randomized Krylov subspace method for $\ell_2-\ell_p$ regularization presented in this paper:
\begin{itemize}
\item Iteratively reweighted norm methods (IRN) which use sketch-to-precondition LSQR to solve the corresponding Tikhonov-regularized subproblems, namely IRN-S2P-LSQR, described in Section \ref{sec:precIRN}.
\item Sketch-and-solve methods based on flexible Krylov subspaces: (hybrid) S\&S-FGMRES and (hybrid) S\&S-FLSQR, or with iterative-reweighted regularization S\&S-IRW-FLSQR and S\&S-IRW-FGMRES, described in Section \ref{sec:flexible_s2s}.
\item Sketch-to-precondition methods based on flexible Krylov subspaces: (hybrid) S2P-FGMRES and (hybrid) S2P-FLSQR, or with iterative-reweighted regularization S2P-IRW-FLSQR and S2P-IRW-FGMRES, described in Section~\ref{sec:flexible_s2p}.
\end{itemize}
Moreover, to evaluate the performance of the new methods, the results obtained using other standard and state-of-the-art algorithms are given for reference:
\begin{itemize}
\item Standard (hybrid) Krylov methods, based on the Arnoldi method (GMRES) and the Golub-Kahan bidiagonalization (LSQR), see e.g. \cite{Chung2024review}, as implemented in \cite{Gazzola2019IRtools}.
\item Iteratively reweighted norm (IRN) methods, where the inner solver for \eqref{eq:reweighted_std} is either based on the Arnoldi factorization (IRN-GMRES), or the Golub-Kahan bidiagonalization  (IRN-LSQR).
\item Variants of flexible Krylov solvers, based on either the flexible Arnoldi or the flexible Golub-Kahan decomposition and, at each iteration, solve \eqref{eq:flex_opt}. In particular, FGMRES and FLSQR, which consider $R_k(\bfy)=0$. The corresponding hybrid versions, which consider $R_k(\bfy)=\|\bfy\|_2^2$: h-FGMRES / h-FLSQR, as described in \cite{Chung2019lp}. And the iteratively-reweighted versions IRW-FGMRES and IRW-FLSQR; which consider $R_k(\bfy)=\|\bfW_k \bfZ_k \bfy\|_2^2$, as described in \cite{Gazzola2021IRW}.
\item Other standard solvers for problems with $\ell_p$ regularization, the fast iterative shrinkage-thresholding algorithm (FISTA) \cite{beck2009fista} and the Sparse Reconstruction by Separable Approximation (SpaRSA) \cite{ Wright2008sparsa}. Both methods require the specification of a regularization parameter a priori, and the implementations that we use are provided in \cite{Gazzola2019IRtools} and by the original authors of \cite{Wright2008sparsa}, respectively.
\end{itemize}

Unless stated otherwise, all the sketching matrices are based on the leverage scores 
of $[\bfA\bfPsi\bar{\bfZ}_k\  \bfb]$ or $[\bfW_k\bar{\bfZ}_k]$,
and the dimension of the sketch is taken to be 4 times the maximum amount of allowed iterations.

The regularization parameters are computed using the DP and the GCV criteria. Moreover, to assess the performance of the proposed methods independently of the regularization parameter choice, we also choose to use the following optimal regularization parameter, i.e.
\begin{equation}\label{eq:opt_reg_p}
    \lambda_k = \arg \min_{\lambda} \|\bfx_k(\lambda)-\bfx_{\rm true}\|_2.
\end{equation}
This is of course not a valid parameter choice in practice, since it requires knowledge about the true solution, but rather a tool for the evaluation of the methods' performance.

The experiments presented in this section were performed using MATLAB on an M-series MacBook Pro. Note that, in this work, we do not present speed comparisons, as the execution time depends heavily on the quality of implementation and hardware.

\subsection{Randomized preconditioning - subset selection example}
\begin{figure}[ht]
    \centering
    \includegraphics[width=\textwidth]{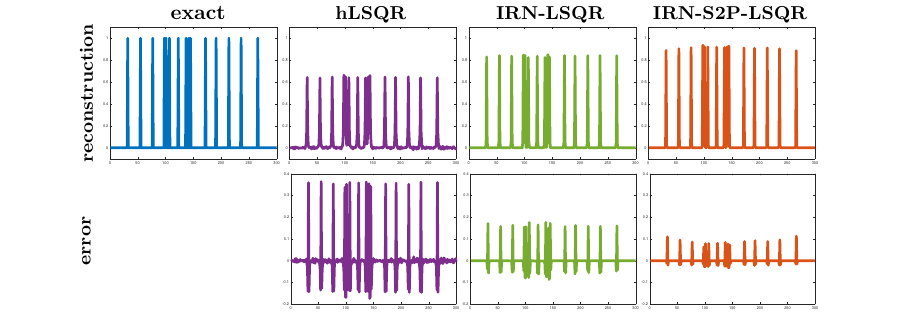}
    \caption{Example 1. Top row: exact and reconstructed coefficients using different methods. Bottom row: errors with respect to the exact solution.}
    \label{ex1:sol}
\end{figure}

This is a simulated example taken from \cite{10.1214/19-STS733} and \cite{10.1214/15-AOS1388} and corresponds to efficient subset selection for linear regression based on $\ell_1$ regularization. The rows of the system matrix $\bfA\in \R^{m \times n}$ are i.i.d. samples from a normal distribution $\mathcal{N}(\bf0,\bfSigma)$ where $\bfSigma \in \R^{n \times n}$ and has entries $\bfSigma_{i,j} =\rho^{|i-j|}$ with $\rho=0.95$ being the predictor correlation level. The entries of the coefficient vector $\bfx_{\rm true}$ (or true solution) are modeled as independent realization of a Bernoulli process with probability parameter $0.1$ (so they are sparse and either 1 or 0). This can be observed in Figure \ref{ex1:sol}, along with the reconstructions computed using different methods and the corresponding errors $\bfx_{\rm true}$$-$$\bfx$. 

\begin{figure}[ht]
\centering
 \centering
    \begin{subfigure}[b]{0.45\textwidth}
        \centering
        \includegraphics[width=4.2cm]{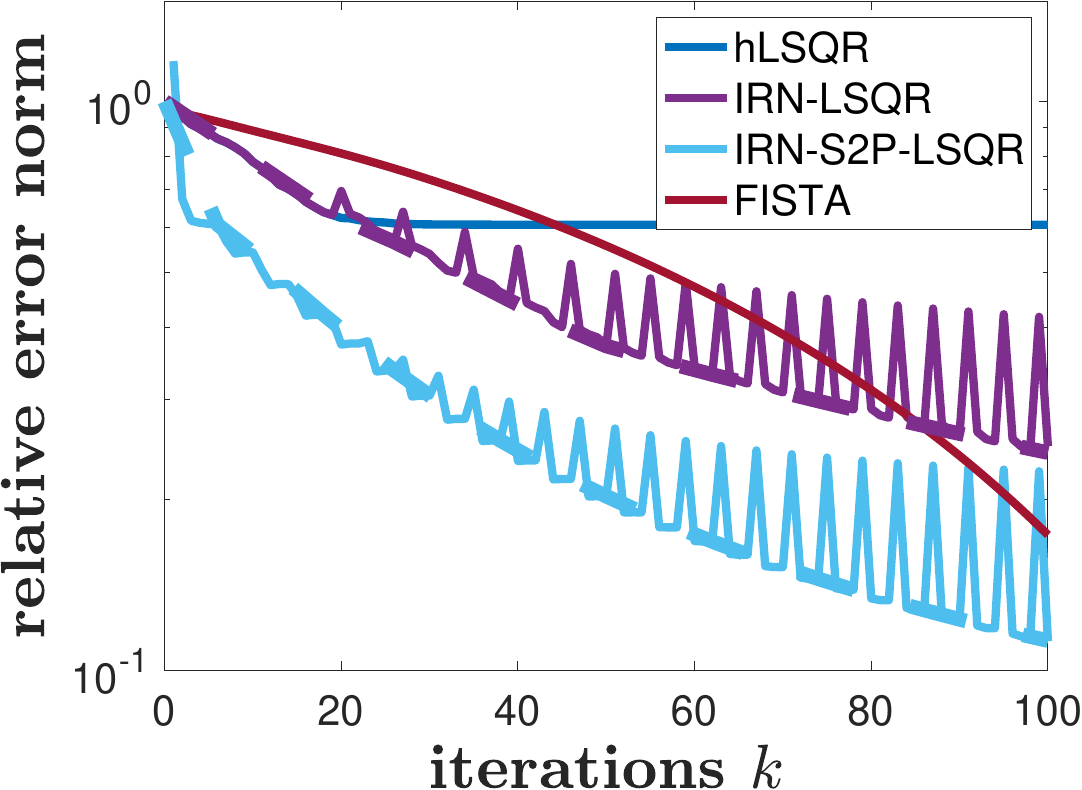}
    \end{subfigure}
    \begin{subfigure}[b]{0.45\textwidth}
        \centering
        \includegraphics[width=4.2cm]{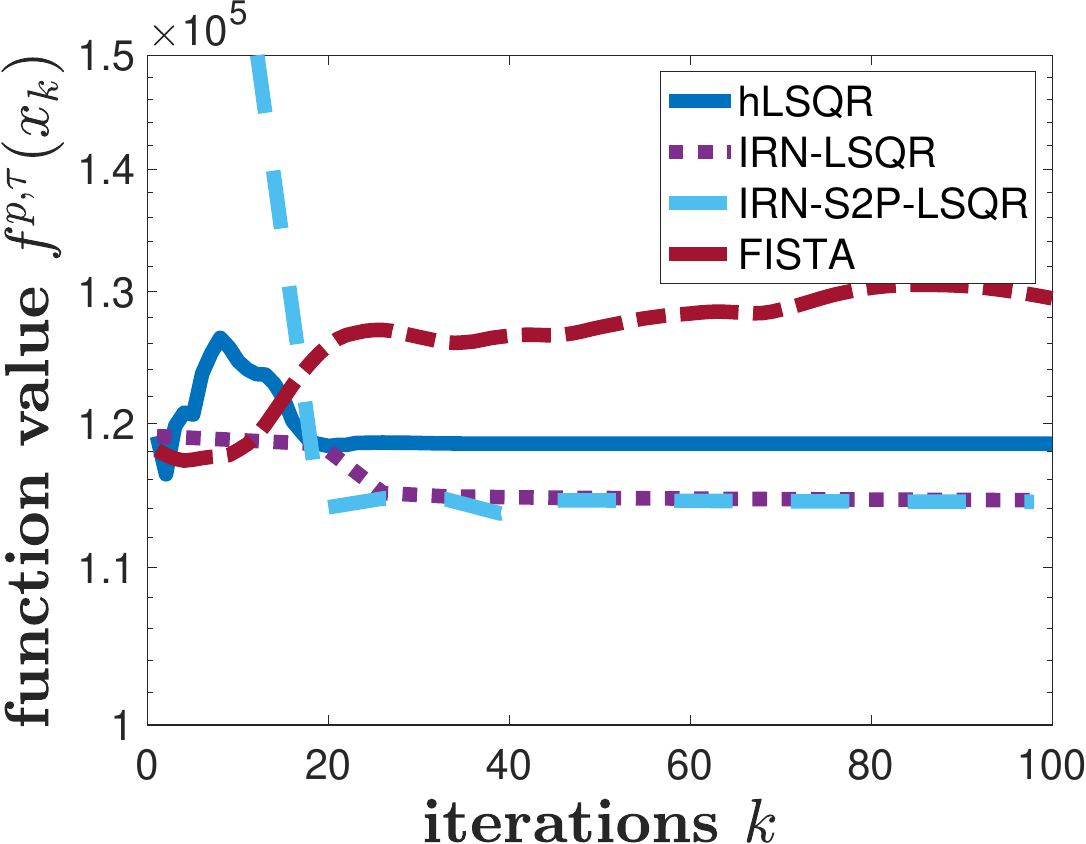}
    \end{subfigure}
    \begin{subfigure}[b]{0.45\textwidth}
        \centering
        \includegraphics[width=4.2cm]{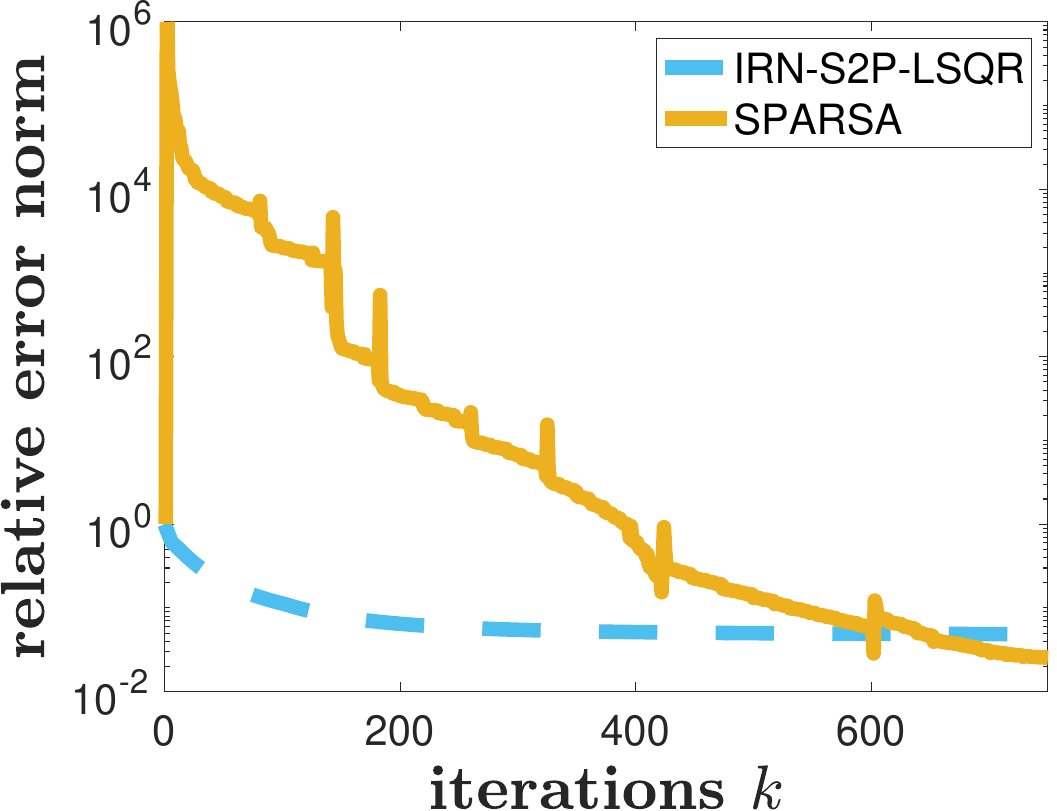}
    \end{subfigure}
    \begin{subfigure}[b]{0.45\textwidth}
        \centering
        \includegraphics[width=4.2cm]{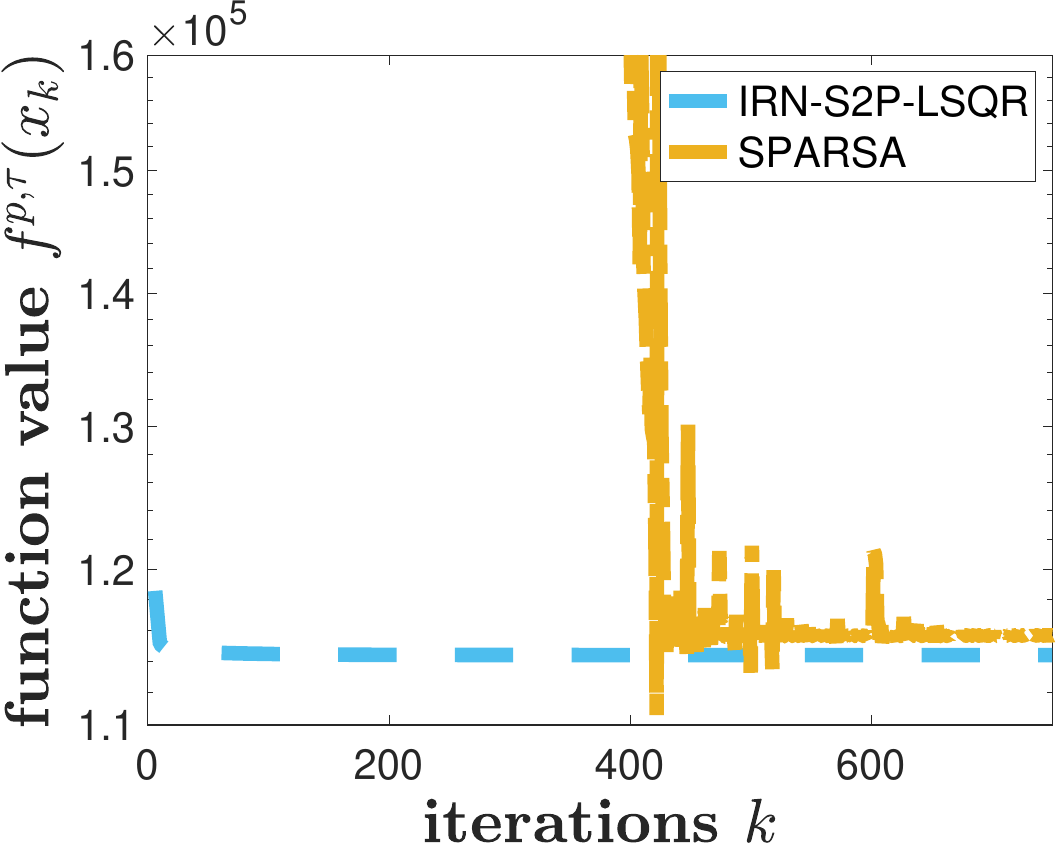}
    \end{subfigure}
    \caption{Example 1. First column: Relative error norm histories for different solvers against the number of iterations (for inner-outer schemes, we display the total number of iterations, not only the outer ones). Second column: value of the function \eqref{eq:smoothed} at each approximate solution for $p=1$ and $\tau=10^{-10}$. Note that on the top left we show the values at each inner-iteration, but for the sake of clarity, the other plots display only the error norm at each outer iteration–indexed by the cumulative count of inner iterations– with linear interpolation applied between successive points.}
    \label{ex1:fixed}
\end{figure}

The error norms for the solutions computed with fixed regularization parameter $\lambda=30000$ can be found in Figure \ref{ex1:fixed} (top row). Note that, in this plot, we show the reconstruction error for all the iterations including the inner iterations of the solvers after the (cold) restarts (after each restart we take $\bfx_0=\bf0$). Note that warm restarts (i.e. after each restart we take $\bfx_0$ to be the previous approximation of the solution) are possible but require technical modifications of the algorithm, see, for example \cite{Chen2014IRLSGS,Onisk2025IRflexible}. However, to simplify the presentation in the following, the plots report the error norm solely at outer iterations, referenced by the total number of inner iterations, while intermediate values are obtained via linear interpolation.

Looking at the error norms in the first column of Figure \ref{ex1:fixed}, one can clearly observe that IRN methods based on Krylov subspaces are very competitive with respect to other solvers, particularly in terms of efficiency. This is particularly true with respect to SpaRSA, which needs 6 times more iterations to obtain a solution of the same quality and therefore is displayed with a different axis, although note that it finally obtains a better solution overall after more iterations. Comparing IRN-LSQR and IRN-S2P-LSQR, i.e. including randomized preconditioning, one can observe that this type of preconditioning really improves the convergence of the method, needing fewer iterations to obtain a comparable solution. This is because LSQR converges faster with IRN-S2P-LSQR, thanks to the high-quality preconditioner.

Moreover, on the second column of Figure \ref{ex1:fixed}, we present the value of the nonlinear regularized least-squares that we are interested in minimizing. One can observe that, besides some small fluctuations (coming from numerical instabilities), this quantity is decreasing monotonically. Moreover, the IRN-S2P-LSQR method achieves the lowest objective function value compared to all the competing methods.

\begin{figure}[ht]
\centering
    \begin{subfigure}[b]{0.45\textwidth}
        \centering
        \includegraphics[height=3.5cm]{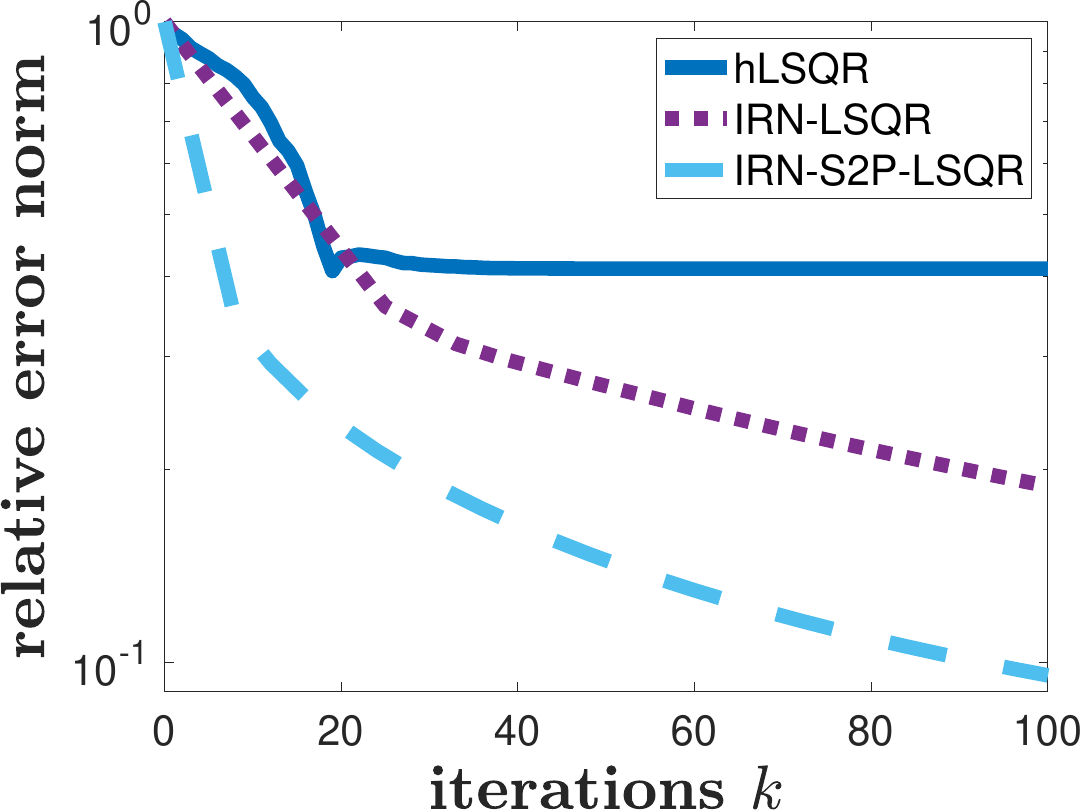}
    \end{subfigure}
    \begin{subfigure}[b]{0.45\textwidth}
        \centering
        \includegraphics[height=3.5cm]{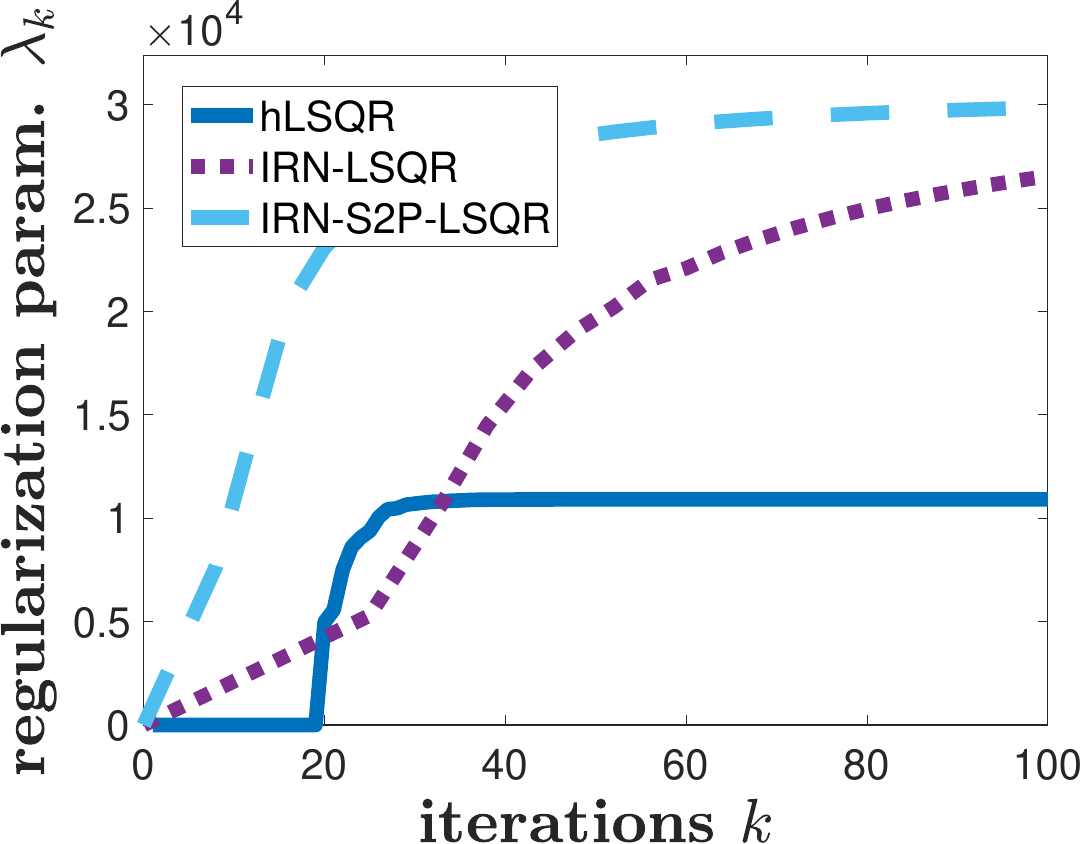}
    \end{subfigure}
    \caption{Example 1. Relative error norm histories (left) and regularization parameters (right) for different solvers against the total number of iterations. The regularization parameters are chosen at each iteration using the discrepancy principle.}
    \label{ex1:dp}
\end{figure}

Last, in Figure \ref{ex1:dp}, we show the relative error norms (left) and regularization parameters (right), for the different methods that allow the regularization parameter to be chosen on-the-fly throughout the iterations. In particular, we use the discrepancy principle in \eqref{eq:dp}. Note that the randomized preconditioning really accelerates the convergence of the IRN method, both in terms of achieving lower relative error norms in less iterations as well as in stabilizing the regularization parameter.  

\subsection{Randomized flexible methods – deblurring example}
\begin{figure}[!ht]
 \centering
    \begin{subfigure}[b]{0.4\textwidth}
        \centering
        \includegraphics[height=3cm]{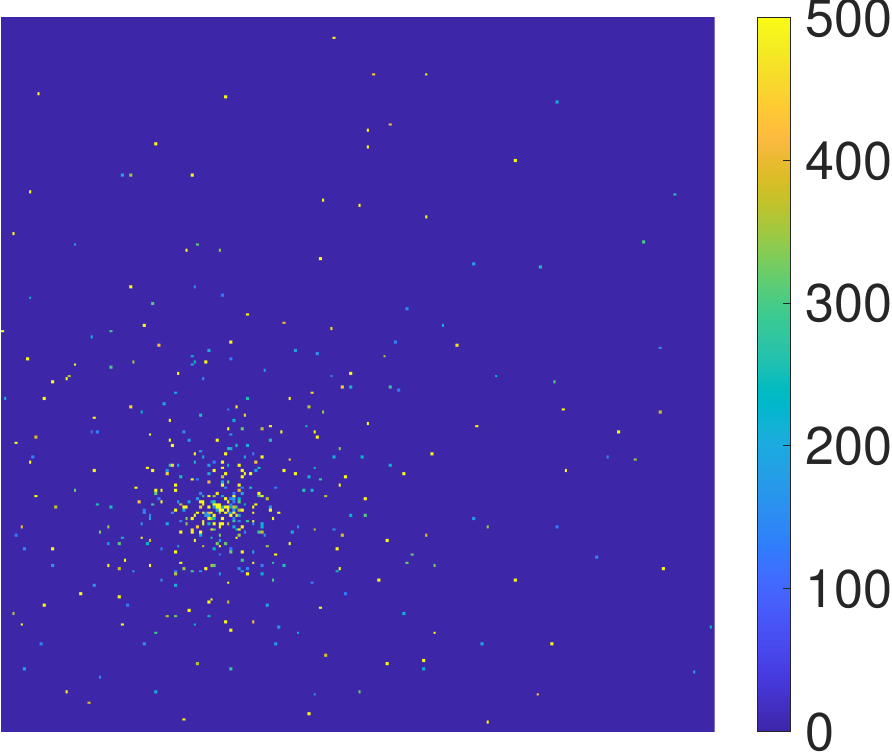}
    \end{subfigure}%
    ~ 
    \begin{subfigure}[b]{0.4\textwidth}
        \centering
        \includegraphics[height=3cm]{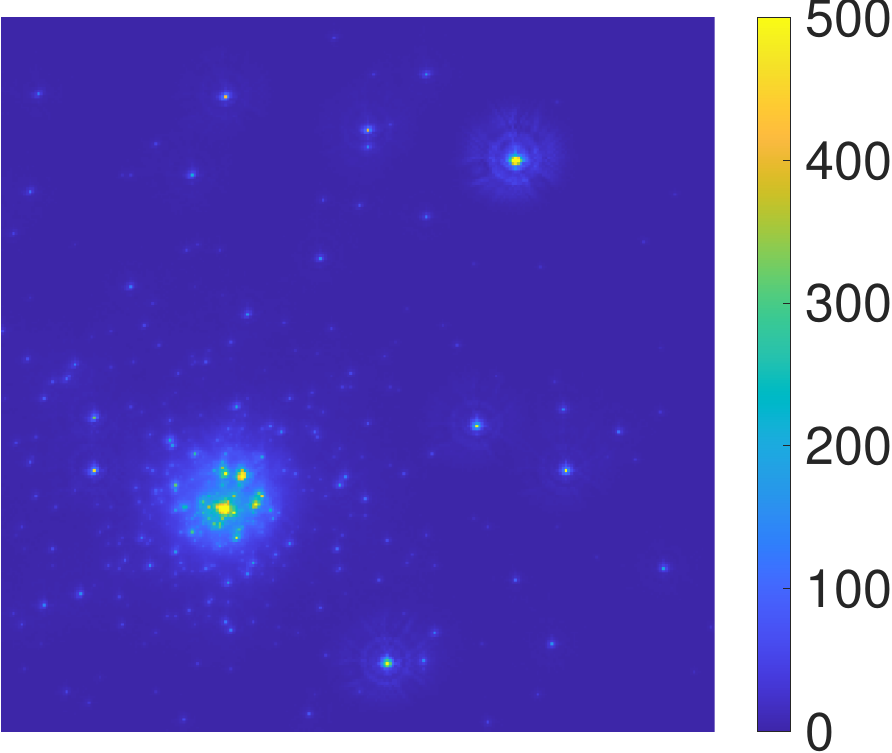}
    \end{subfigure}
    \caption{Example 2. Left: Exact image. Right: Noisy, blurred image.}
    \label{fig:ex3_setting}
\end{figure}

This numerical experiment corresponds to a deblurring problem simulating the imaging of stars under atmospheric blur, where the solution has $256\times 256$ pixels (so that $\bfx^{65536}$). Even if the corresponding (square) matrix models a non-homogeneous blur, it is highly structured, so one can perform matrix-vector products with $\bfA$ and its transpose efficiently. For more details, see \cite{Nagy1998degraded}. The solution of the noiseless system, which can be found in  \cite{RestoreTools}, is sparse by construction in the sense that $92.8 \%$ of the pixels have a magnitude below $10^{-10}$, and can be observed in Figure \ref{fig:ex3_setting}. In the same figure, one can also observe the measurements, which have been corrupted with Gaussian white noise with noise level $0.01$.

\begin{figure}[!ht]
 \centering
  \includegraphics[width=\textwidth]{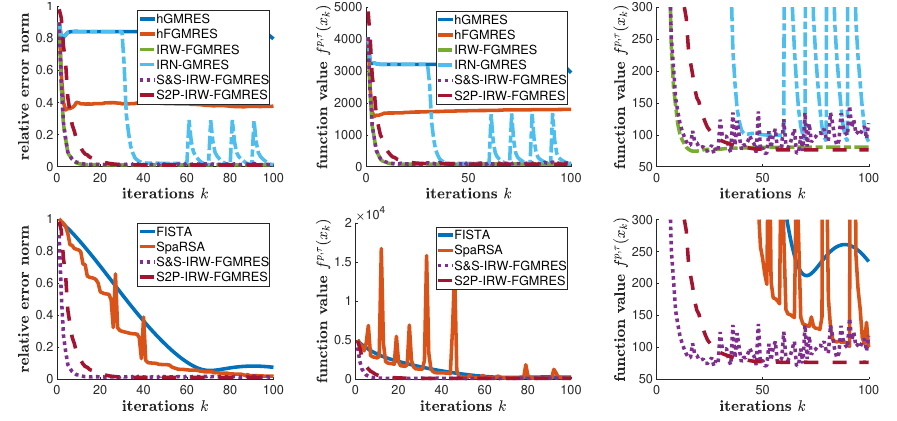}
    \caption{Example 2. First column: Relative error norm histories for randomized flexible iteratively reweighted method. Second column: Value of the non-linear objective function \eqref{eq:smoothed} with $p=1$ and $\tau=10^{-10}$ at the approximate solutions obtained using different methods. Third column: close-up of the bottom of the plots in the second column. The new methods are compared to other standard and flexible Krylov solvers (top row), as well as other standard solvers for $\ell_1$ regularization (bottom row).}
    \label{fig:Ex3_Enrm_Fx}
\end{figure}

Since the system matrix $\bfA$ is square, we use this example to test and compare methods based on the (possibly flexible) Arnoldi method. Due to the level of ill-posedness of this example, we only compare methods with explicit regularization. In particular, we fix the regularization parameter $\lambda = 0.1855$, which we find to be a good choice for the  original full dimensional problem. The first column of Figure~\ref{fig:Ex3_Enrm_Fx} displays the relative error norm histories of the reconstructed solutions obtained using the new methods (based on a truncated flexible Arnoldi process with truncation parameter $\ell=4$) compared to other (flexible) Krylov solvers (top) and other methods tailored to problems with $\ell_1$ regularization (bottom). One can clearly observe that, for this example, using a sketch-and-solve approach in combination with the iteratively reweighted flexible GMRES (S\&S-IRW-FGMRES) shows the fastest convergence, along with the most expensive IRW-FGMRES, which solves the projected problems exactly. These are closely followed by using sketch-to-precondition. It is also evident to see that methods that consider an explicit standard Tikhonov regularization (hGMRES and hFGMRES) instead of approximating an $\ell_1$, drastically under-perform for this experiment.

On the second column (and third for a close-up) of Figure \ref{fig:Ex3_Enrm_Fx}, one can observe the values of the original non-linear objective function \eqref{eq:smoothed} with $p=1$ and $\tau=10^{-10}$ at every iteration. One can see that, unsurprisingly, the sketch-to-precondition approach achieved the lowest values overall, and it is almost monotonic in all iterations (with small perturbations induced by computational errors) as expected from Proposition~2. Note that this is different from the sketch-and-solve approach.

\begin{figure}[!ht]
 \centering
    \begin{subfigure}[b]{0.5\textwidth}
        \centering
        \includegraphics[height=3.5cm]{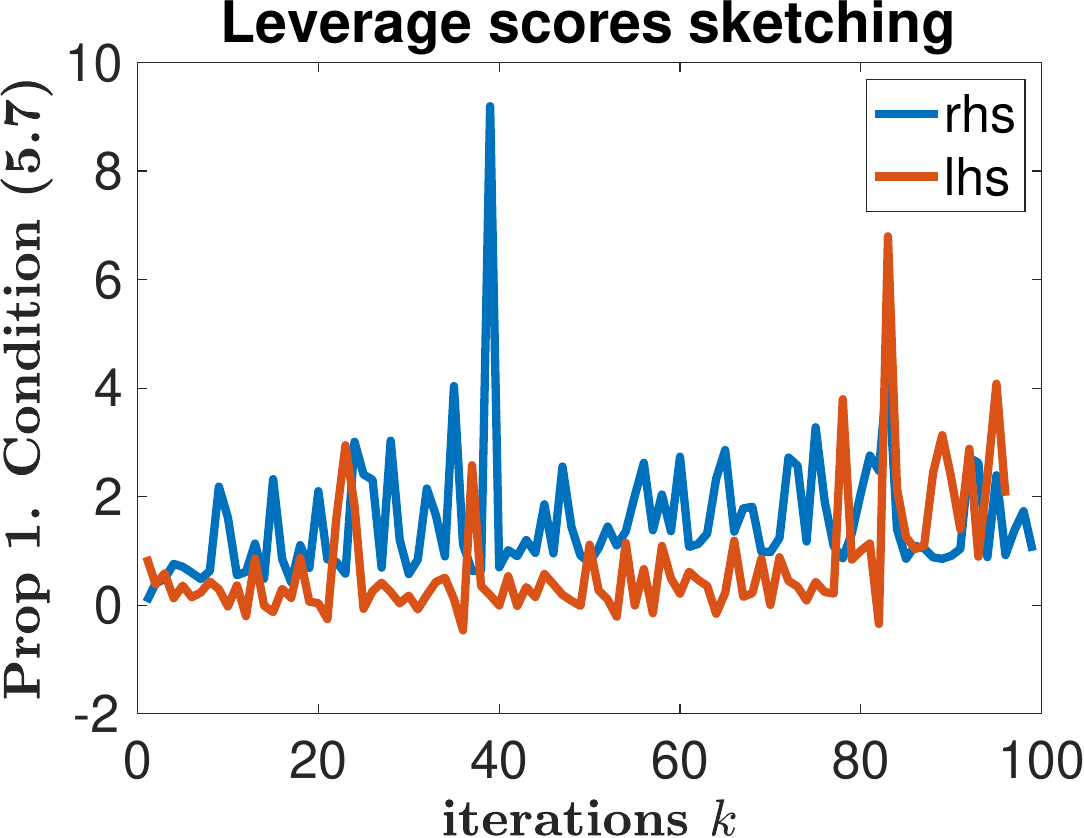}
    \end{subfigure}%
    \begin{subfigure}[b]{0.5\textwidth}
        \centering
        \includegraphics[height=3.5cm]{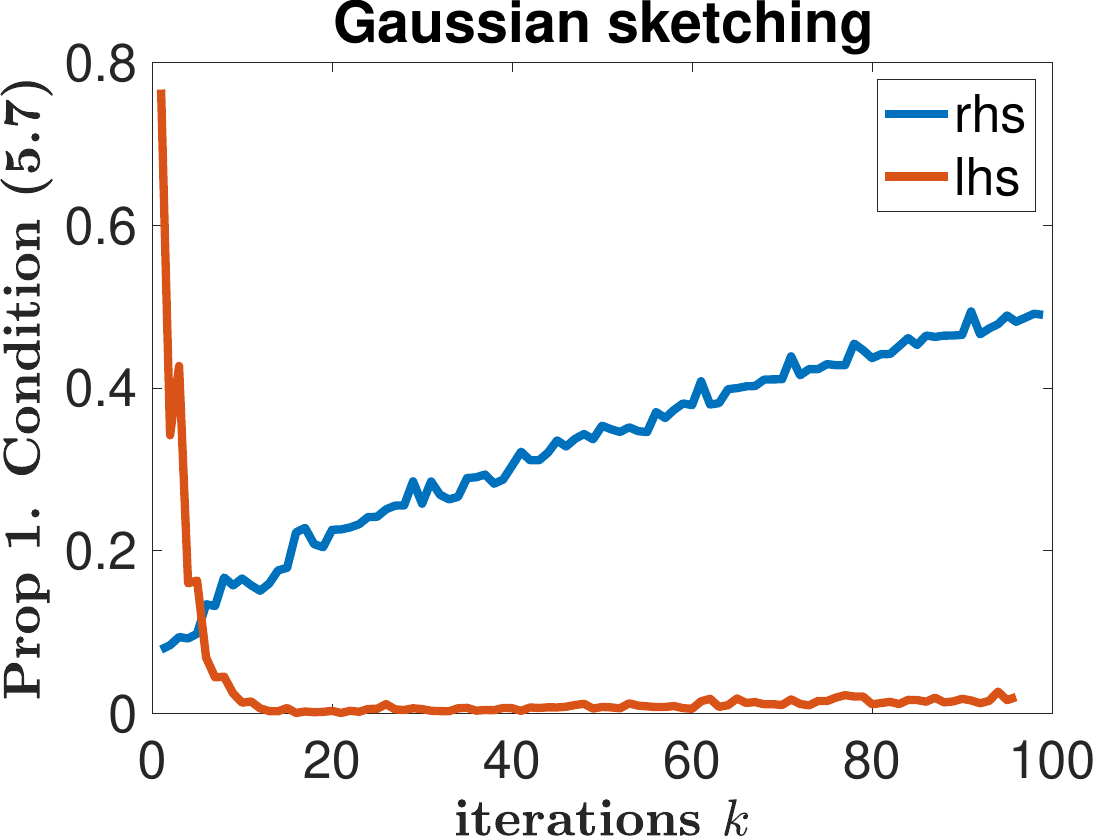}
    \end{subfigure}
   \caption{Example 2. Proposition~1 empirically for two types of sketch.}
    \label{fig:ex2_lemma1}
\end{figure}

Last, to study the convergence of the sketch-and-solve approach, in Figure \ref{fig:ex2_lemma1} we show some empirical results illustrating the behaviour of the condition for monotonicity \eqref{eq:mono_cond} stated Proposition~1. As one can observe, this condition is only met in the first few iterations of the algorithm. First, it is important to mention that this  condition only covers for the worst-case-scenario in the sketching error bounds. However, it is a topic of future research to use this information to devise provably converging sketch-and-solve methods based on restarts, see \cite{Onisk2025IRflexible}, sketch refinement, or iterative sketching \cite{Epperly2024Stable}.

\subsection{Randomized flexible methods – tomography example}
\begin{figure}[!ht]
 \centering
    \begin{subfigure}[b]{0.4\textwidth}
        \centering
        \includegraphics[height=3cm]{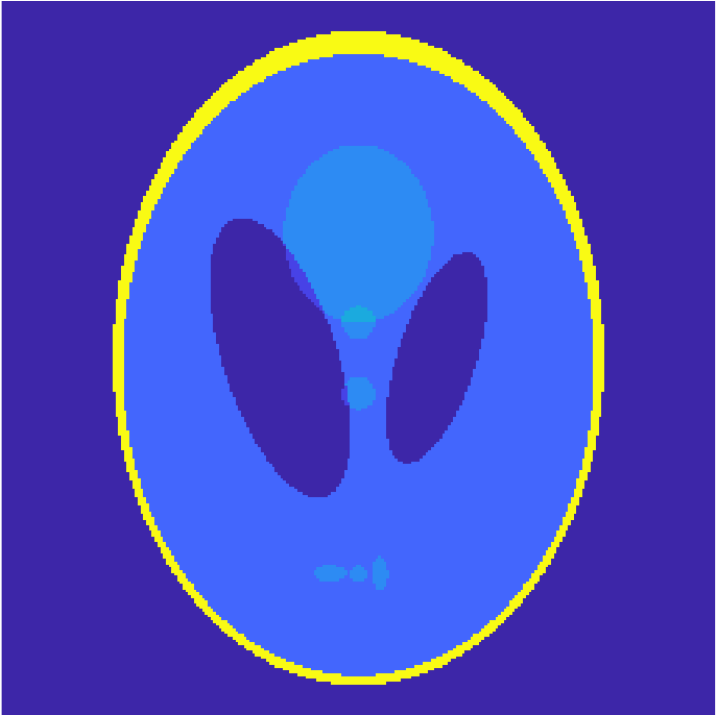}
    \end{subfigure}%
    ~ 
    \begin{subfigure}[b]{0.4\textwidth}
        \centering
        \includegraphics[height=3cm]{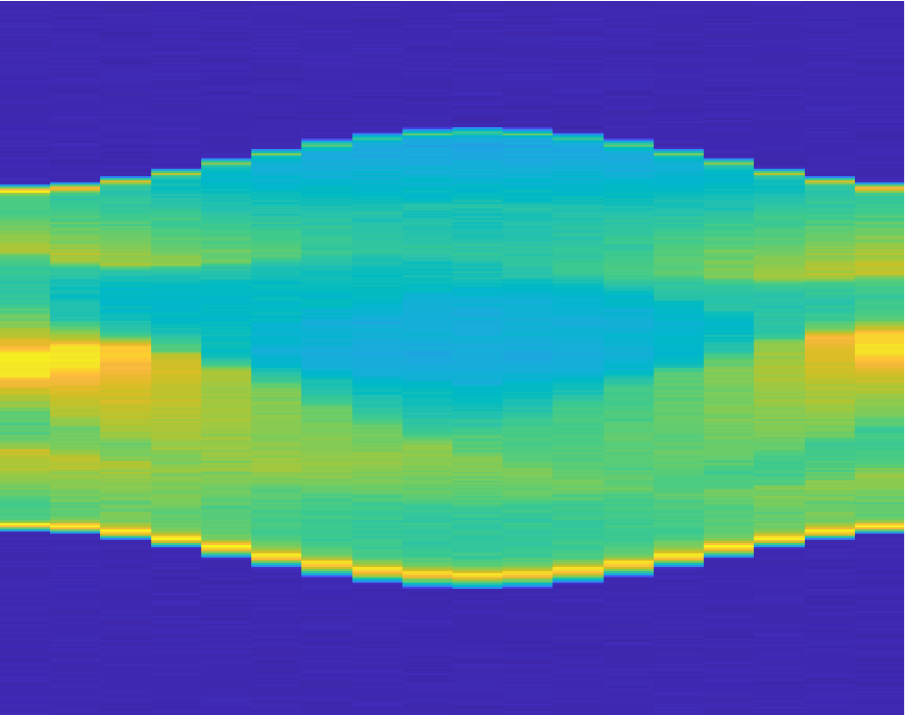}
    \end{subfigure}
    \caption{Example 3. Left: Exact image. Right: Noisy sinogram.}
    \label{fig:ex2_setting}
\end{figure}

This experiment consists of a test CT problem with undersampled projection angles concerning the Shepp Logan phantom (assumed to span $256 \times 256$ pixels). More specifically, this synthetic example assumes a parallel geometry with 362 rays,  18 equispaced projection angles between 1 and 180 degrees, and it is generated using IRtools \cite{Gazzola2019IRtools}. The corresponding  underdetermined system is then $\bfA \in \mathbb{R}^{6516 \times 65536}$. Gaussian white noise with a noise level of 0.01 is added to the measurements. The exact solution and the noisy measurements (a.k.a. sinogram in the CT context), can be observed in Figure \ref{fig:ex2_setting}.

\begin{figure}[!ht]
    \centering
    \includegraphics[height=4.5cm]{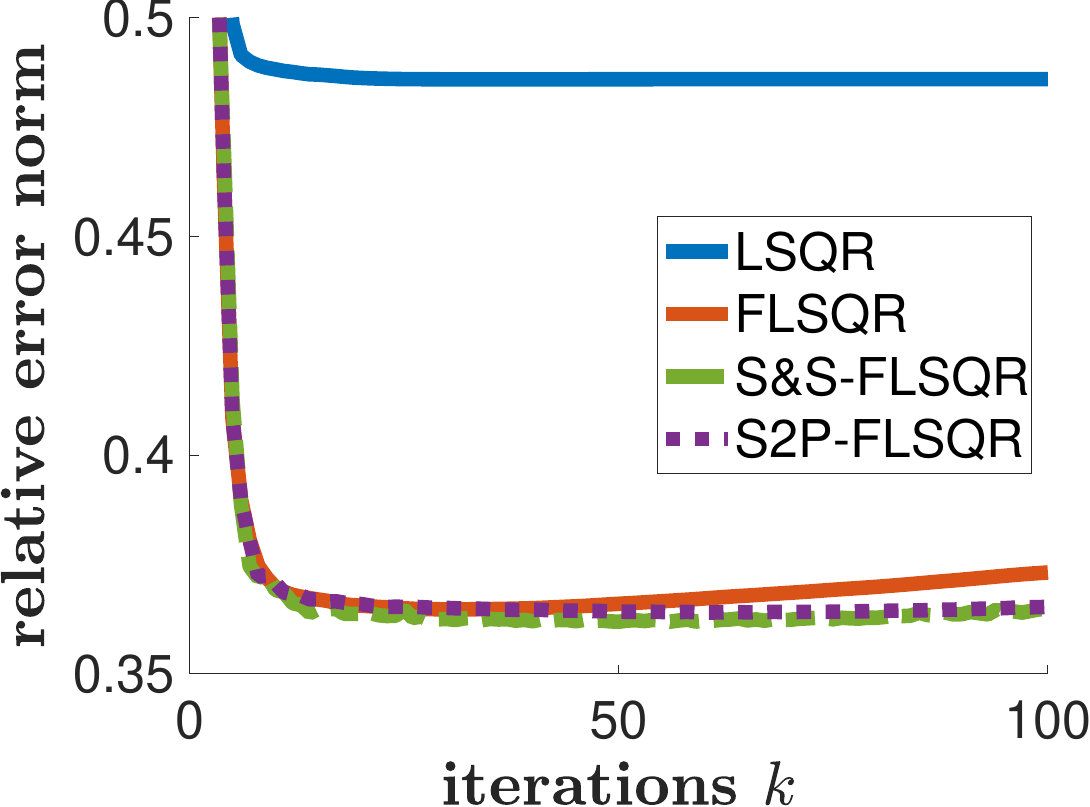}
    \caption{Example 3. Relative error norm histories for different (potentially randomized and/or flexible) Krylov solvers without explicit regularization.}
    \label{fig:Enrm_no_reg}
\end{figure}

\begin{figure}[!ht]
    \centering
    \includegraphics[width=\textwidth]{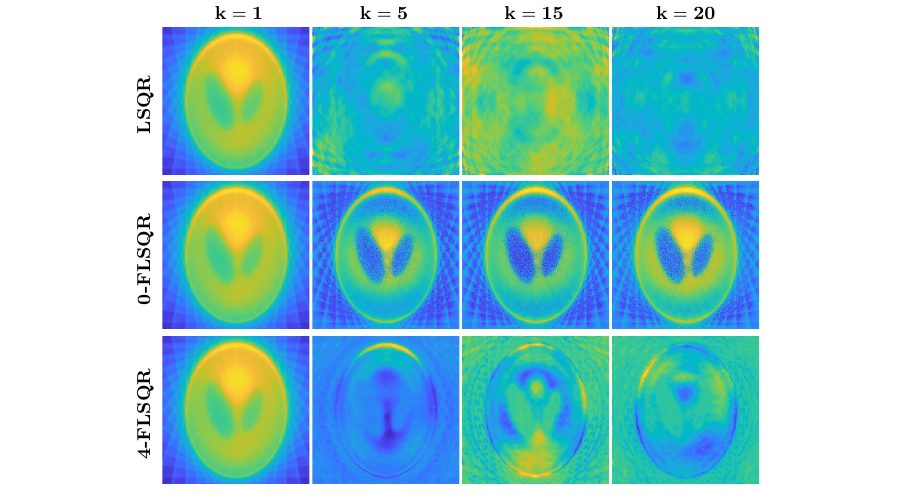}
    \caption{Example 3. Different basis vectors corresponding to the $k^{th}$ columns of the appropriate matrix constructed using Arnoldi (upper row), $0$-truncated flexible Arnoldi (middle row) and $4$-truncated flexible Arnoldi (lower row).}
    \label{fig:basis_vectors}
\end{figure}

For this experiment, we compare first the relative error histories for different methods that do not consider explicit regularization (i.e., where $\lambda=0$), namely LSQR, flexible LSQR (FLSQR), sketch-and-solve flexible LSQR (S\&S-FLSQR) and sketch-to-precondition flexible LSQR  (S2P-FLSQR). Note that the two latter examples are based on a truncated flexible GK process with truncation parameter $\ell=4$. Note that, even if these methods do not have explicit regularization, implicit regularization can be achieved in the flexible methods by means of the preconditioning. The effect of using a better space for the solution can be dramatic, as one can see in Figure \ref{fig:Enrm_no_reg} comparing methods based on the classic Golub-Kahan bidiagonalization and its flexible Golub-Kahan counterpart. One can observe the basis vectors corresponding to different methods in Figure \ref{fig:basis_vectors}. First, note that the GKB process used in LSQR produces highly oscillatory basis vectors, whereas the vectors corresponding to new (truncated) flexible methods incorporate more information about the image `sparsity' patterns. Also note that the truncation in the orthogonalization of the basis has a great impact in the basis vectors' appearance. 

\begin{figure}[!ht]
 \centering     
    \includegraphics[width=\textwidth]{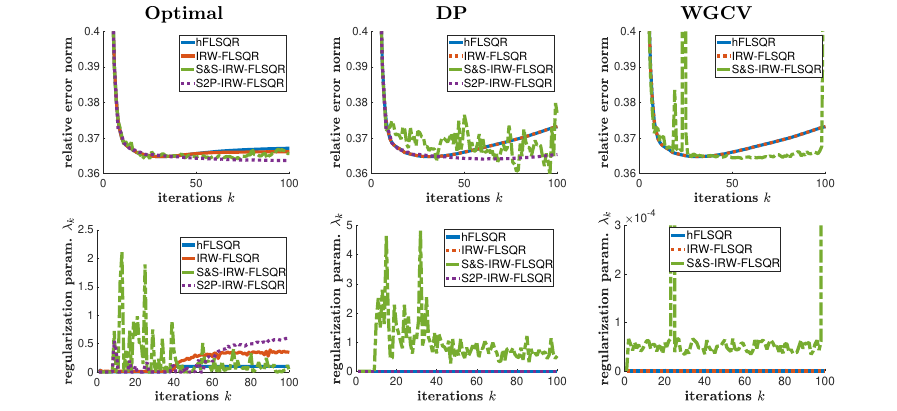}
    \caption{Example 3. Top row: Relative error norm histories for different methods using explicit regularization. Bottom row: Regularization parameters chosen at each iteration. Each column corresponds to a different regularization parameter choice criterion: optimal regularization parameter according to \eqref{eq:opt_reg_p}, discrepancy principle (DP) and weighted generalized cross validation (WGCV), respectively.}
    \label{fig:ex2_reg}
\end{figure}

To compare the performance of the different methods including explicit regularization, the first column of Figure~\ref{fig:ex2_reg} shows the relative error norm histories (top) and the corresponding optimal regularization parameters chosen at each iteration according to \eqref{eq:opt_reg_p} (bottom). One can observe that the sketched methods have a very comparable performance to the flexible non-sketched counterparts. In particular, these methods seem to perform slightly better for this example, due to the differences in the basis vectors. Here, however, it is important to stress that one of the main advantages of the sketched methods is related to a higher computational efficiency. Moreover, to evaluate the effectiveness of different regularization parameter choices, the second and third columns of Figure~\ref{fig:ex2_reg} show the relative error norms for two different regularization parameter choice criteria.  In particular, we compare the discrepancy principle (DP) and the weighted generalized cross-validation (WGCV), in the second and third columns respectively. First, note that the discrepancy principle is not as effective for the sketch-and-solve flexible method. This is not surprising, since even if we can observe that the approximated solution using sketching is close to the approximated solution without sketching, the sketched residual norm might still not a good enough approximation to the true residual norm. Also note that the WGCV is not used in the sketch-to-precondition method, since this requires estimating a trace, which is not within the scope of this paper.

\section{Conclusions}
In this paper we have presented an overview of how randomized numerical linear algebra can be used in combination with (potentially flexible) Krylov subspace methods to solve linear discrete inverse problems with $\ell_1$ regularization.

All the methods presented here are based on a majorization-minimization scheme of (a smoothed version) of the original functional, where a sequence of least-squares problems with (weighted) Tikhonov regularization needs to be (approximately) solved at each iteration.

On one hand, when the system matrix is tall and skinny (or short and fat), one can directly use effective preconditioning based on sketching to accelerate the convergence of the methods. Moreover, in the numerical experiments, we also observe that this technique might help overcome stagnation of the non-preconditioned methods.

On the other hand, for closer to square systems, methods based on flexible Krylov subspaces can be used in combination with sketch-and-solve or sketch-to-precondition in the projected subproblems. In this paper, we analyse the computational cost of these methods, and we show that the generated subspaces using truncated flexible Krylov solvers are not equivalent to their non-truncated counterparts (differently from standard Krylov solvers). Moreover, we give theoretical insight (and conditions) for the monotonicity of the original non-linear objective function with respect to the iterations in each case. 

In the numerical experiments, it is evident that sketching renders traditional methods more computationally efficient at a minimal cost of accuracy. More importantly, sketching does not seem to amplify the noise in linear discrete inverse problems, becoming a very viable alternative to other solvers for problems with sparse solutions. 

\bibliographystyle{abbrv}
\bibliography{biblio}

\begin{thebibliography}{10}

\bibitem{beck2009fista}
A.~Beck and M.~Teboulle.
\newblock A fast iterative shrinkage-thresholding algorithm for linear inverse
  problems.
\newblock {\em SIAM J. Imaging Sci.}, 2(1):183--202, 2009.

\bibitem{Bjorck1996LS}
{\AA}.~Bj{\" o}rck.
\newblock {\em Numerical methods for least squares problems}.
\newblock SIAM, Philadelphia, PA, 1996.

\bibitem{Callaghan2008IllPosed}
K.~J. Callaghan and J.~Chen.
\newblock Revisiting the collinear data problem: An assessment of estimator
  ``ill-conditioning'' in linear regression.
\newblock {\em Pract. Assess. Res. Eval.}, 13, 1 2008.

\bibitem{calvetti2005priorconditioners}
D.~Calvetti and E.~Somersalo.
\newblock Priorconditioners for linear systems.
\newblock {\em Inverse problems}, 21(4):1397--1418, 2005.

\bibitem{Candes2007Sparsity}
E.~Candès, M.~Wakin, and S.~Boyd.
\newblock Enhancing sparsity by reweighted {L1} minimization.
\newblock {\em J. Fourier Anal. Appl.}, 14:877--905, 11 2007.

\bibitem{Chen2014IRLSGS}
C.~Chen, J.~Huang, L.~He, and H.~Li.
\newblock Fast iteratively reweighted least squares algorithms for
  analysis-based sparse reconstruction.
\newblock {\em Med. Image Anal.}, 49, 11 2014.

\bibitem{Chung2019lp}
J.~Chung and S.~Gazzola.
\newblock Flexible {K}rylov methods for $\ell_p$ regularization.
\newblock {\em SIAM J. Sci. Comput.}, 41(5):S149--S171, 2019.

\bibitem{Chung2024review}
J.~Chung and S.~Gazzola.
\newblock Computational methods for large-scale inverse problems: A survey on
  hybrid projection methods.
\newblock {\em SIAM Rev.}, 66(2):205--284, 2024.

\bibitem{chung2025randomized}
J.~Chung and S.~Gazzola.
\newblock Randomized {K}rylov methods for inverse problems, 2025.

\bibitem{Daubechies2010IRLS}
I.~Daubechies, R.~DeVore, M.~Fornasier, and C.~S. Güntürk.
\newblock Iteratively reweighted least squares minimization for sparse
  recovery.
\newblock {\em Comm. Pure Appl. Math.}, 63(1):1--38, 2010.

\bibitem{10.1214/15-AOS1388}
A.~K. Dimitris~Bertsimas and R.~Mazumder.
\newblock {Best subset selection via a modern optimization lens}.
\newblock {\em Ann. Statist.}, 44(2):813--852, 2016.

\bibitem{drineas2012fast}
P.~Drineas, M.~Magdon-Ismail, M.~W. Mahoney, and D.~P. Woodruff.
\newblock Fast approximation of matrix coherence and statistical leverage.
\newblock {\em J. Mach. Learn. Res.}, 13(1):3475--3506, 2012.

\bibitem{Epperly2024Stable}
E.~N. Epperly.
\newblock Fast and forward stable randomized algorithms for linear
  least-squares problems.
\newblock {\em SIAM J. Matrix Anal. Appl.}, 45(4):1782--1804, 2024.

\bibitem{Gazzola2019IRtools}
S.~Gazzola, P.~C. Hansen, and J.~G. Nagy.
\newblock {IR T}ools: a {MATLAB} package of iterative regularization methods
  and large-scale test problems.
\newblock {\em Numer. Algorithms}, 81(3):773–811, July 2019.

\bibitem{Gazzola2014FAT}
S.~Gazzola and J.~Nagy.
\newblock Generalized {A}rnoldi-{T}ikhonov method for sparse reconstruction.
\newblock {\em SIAM J. Sci. Comput.}, 36, 01 2014.

\bibitem{Gazzola2021IRW}
S.~Gazzola, J.~G. Nagy, and M.~S. Landman.
\newblock Iteratively reweighted {FGMRES} and {FLSQR} for sparse
  reconstruction.
\newblock {\em SIAM J. Sci. Comput.}, 43(5):S47--S69, 2021.

\bibitem{Sabate2019TV}
S.~Gazzola and M.~{Sabaté Landman}.
\newblock Flexible {GMRES} for total variation regularization.
\newblock {\em BIT}, 59(3):721--746, Sept. 2019.

\bibitem{Gazzola2020param}
S.~Gazzola and M.~Sabaté~Landman.
\newblock Krylov methods for inverse problems: Surveying classical, and
  introducing new, algorithmic approaches.
\newblock {\em GAMM-Mitt.}, 43(4):e202000017, 2020.

\bibitem{jimaging7100216}
S.~Gazzola, S.~J. Scott, and A.~Spence.
\newblock Flexible {K}rylov methods for edge enhancement in imaging.
\newblock {\em J. Imaging}, 7(10), 2021.

\bibitem{Hansen2010}
P.~C. Hansen.
\newblock {\em Discrete Inverse Problems: Insight and Algorithms}.
\newblock SIAM, Philadelphia, 2010.

\bibitem{10.1214/19-STS733}
T.~Hastie, R.~Tibshirani, and R.~Tibshirani.
\newblock {Best Subset, Forward Stepwise or Lasso? Analysis and Recommendations
  Based on Extensive Comparisons}.
\newblock {\em Statist. Sci.}, 35(4):579 -- 592, 2020.

\bibitem{Huang2017MM}
G.~Huang, A.~Lanza, S.~Morigi, L.~Reichel, and F.~Sgallari.
\newblock Majorization–minimization generalized {K}rylov subspace methods for
  $\ell_p$ – $\ell_q$ optimization applied to image restoration.
\newblock {\em BIT}, 57, 01 2017.

\bibitem{mm}
D.~R. Hunter and K.~Lange.
\newblock A tutorial on {MM} algorithms.
\newblock {\em Amer. Statist.}, 58(1):30--37, 2004.

\bibitem{kutner2004}
M.~Kutner, C.~Nachtsheim, and J.~Neter.
\newblock {\em Applied Linear Regression Models}.
\newblock Irwin/McGraw-Hill series in operations and decision sciences.
  McGraw-Hill/Irwin, 2004.

\bibitem{Lanza2015GKS}
A.~Lanza, S.~Morigi, L.~Reichel, and F.~Sgallari.
\newblock A generalized {K}rylov subspace method for $\ell_p$-$\ell_q$
  minimization.
\newblock {\em SIAM J. Sci. Comput.}, 37(5):S30--S50, 2015.

\bibitem{Kolda2022sampling}
B.~W. Larsen and T.~G. Kolda.
\newblock Sketching matrix least squares via leverage scores estimates, 2022.

\bibitem{Martinsson_Tropp_2020}
P.-G. Martinsson and J.~A. Tropp.
\newblock Randomized numerical linear algebra: Foundations and algorithms.
\newblock {\em Acta Numer.}, 29:403–572, 2020.

\bibitem{meier2022randomized}
M.~Meier and Y.~Nakatsukasa.
\newblock Randomized algorithms for {T}ikhonov regularization in linear least
  squares.
\newblock {\em arXiv preprint arXiv:2203.07329}, 2022.

\bibitem{Morozov1966OnTS}
V.~A. Morozov.
\newblock On the solution of functional equations by the method of
  regularization.
\newblock {\em Dokl. Math.}, 7:414--417, 1966.

\bibitem{Nagy1998degraded}
J.~G. Nagy and D.~P. O'{L}eary.
\newblock Restoring images degraded by spatially-variant blur.
\newblock {\em SIAM J. Sci. Comput.}, 19:1063--1082, 1998.

\bibitem{RestoreTools}
J.~G. Nagy, K.~Palmer, and L.~Perrone.
\newblock Iterative methods for image deblurring: {A} {MATLAB} object-oriented
  approach.
\newblock {\em Numer. Algorithms}, 36:73–93, 2004.

\bibitem{Onisk2025IRflexible}
L.~Onisk and M.~S. Landman.
\newblock Iterative refinement and flexible iteratively reweighed solvers for
  linear inverse problems with sparse solutions, 2025.

\bibitem{renaut2017hybrid}
R.~A. Renaut, S.~Vatankhah, and V.~E. Ardestani.
\newblock Hybrid and iteratively reweighted regularization by unbiased
  predictive risk and weighted {GCV} for projected systems.
\newblock {\em SIAM J. Sci. Comput.}, 39(2):B221--B243, 2017.

\bibitem{rodriguez2008IRN}
P.~Rodr\'{i}guez and B.~Wohlberg.
\newblock An efficient algorithm for sparse representations with $\ell^{p}$
  data fidelity term.
\newblock In {\em Proceedings of 4th IEEE Andean Technical Conference
  (ANDESCON)}, Cusco, Per\'{u}, Oct. 2008.

\bibitem{Sabate2025randCMRH}
M.~Sabat{\'e}~Landman, A.~N. Brown, J.~Chung, and J.~G. Nagy.
\newblock Randomized and inner-product free {K}rylov methods for large-scale
  inverse problems.
\newblock {\em Numer. Algorithms}, 2025.

\bibitem{Sabate2024gmd}
M.~Sabat\'e~Landman, J.~Chung, J.~Jiang, S.~M. Miller, and A.~K. Saibaba.
\newblock A joint reconstruction and model selection approach for large-scale
  linear inverse modeling (mshybr v2).
\newblock {\em Geosci. Model Dev.}, 17(23):8853--8872, 2024.

\bibitem{Tibshirani1996Lasso}
R.~Tibshirani.
\newblock Regression shrinkage and selection via the lasso.
\newblock {\em J. R. Stat. Soc. B.}, 58(1):267--288, 1996.

\bibitem{Wright2008sparsa}
S.~J. Wright, R.~D. Nowak, and M.~A.~T. Figueiredo.
\newblock Sparse reconstruction by separable approximation.
\newblock In {\em 2008 IEEE International Conference on Acoustics, Speech and
  Signal Processing}, pages 3373--3376, 2008.

\end{thebibliography}
\end{document}